\documentclass[%
 aip,
 sd,%
 amsmath,amssymb,
 reprint,%
]{revtex4-1}

\usepackage{graphicx}% Include figure files
\usepackage{dcolumn}% Align table columns on decimal point
\usepackage{bm}% bold math
\usepackage{graphicx}% Include figure files
\usepackage{dcolumn}% Align table columns on decimal point
\usepackage{bm}% bold math
\usepackage{color}
\usepackage{hyperref}% add hypertext capabilities
\usepackage[mathlines]{lineno}% Enable numbering of text and display math
\newcommand{\RNum}[1]{\uppercase\expandafter{\romannumeral #1\relax}}

\def\half{\frac{1}{2}}

\def\dsp{\displaystyle}
\def\cal{\mathcal}
\def\bm{\boldsymbol}
\def\wh{\widehat}
\def\wt{\widetilde}
\def\eg{{\it e.g.}}

\def\red{ }
\def\blue{ }

\newtheorem{Thm}{Theorem}

\newtheorem{rem}[Thm]{Remark}

\begin{document}

\title{Computation of the Memory Functions in the Generalized Langevin Models
for Collective Dynamics of Macromolecules}% Force line breaks with \\
%\thanks{Footnote to title of article.}
\author{Minxin Chen}
 \email{chenmx@gmail.com}
\affiliation{%
Center for System Biology,
Department of Mathematics, Soochow University, Suzhou 215006,
China %\\This line break forced% with \\
}%
\author{Xiantao Li}
\email{xli@math.psu.edu}%Lines break automatically or can be forced with \\
\author{Chun Liu}%
 \email{liu@math.psu.edu}
\affiliation{Department of Mathematics,
Pennsylvania State University, University Park, PA 16802, US %\\This line break forced with \textbackslash\textbackslash
}%

%\author{Minxin Chen}
% \affiliation{Center for System Biology,
%Department of Mathematics, Soochow University, Suzhou 215006,
%China }
%\email{chenmx@gmail.com}
%
%\author{Xiantao Li}
% \affiliation{Department of Mathematics,
%Pennsylvania State University, University Park, PA 16802, US }
%\email{xli@math.psu.edu, liu@math.psu.edu}
%
%\author{Chun Liu}
%\affiliation{%
%Department of Mathematics,
%Pennsylvania State University, University Park, PA 16802, US %\\This line break forced% with \\
%}%

\date{\today}% It is always \today, today,

\begin{abstract}
We present a numerical method to compute the approximation of the memory functions in the generalized Langevin
models for collective dynamics of macromolecules. We first derive the exact expressions of the memory
functions, obtained from projection to subspaces that correspond to the selection of coarse-grain variables.
In particular, the memory functions are expressed in the forms of  matrix functions, which will then be approximated by Krylov-subspace
methods. It will also be demonstrated that the random noise can be approximated under the same framework, and the
fluctuation-dissipation theorem is automatically satisfied.
The accuracy of the method is examined through several numerical examples.
\end{abstract}

\pacs{{\red 05.40.Ca, 05.70.-a, 05.10.G}}% PACS, the Physics and Astronomy
\keywords{Generalized Langevin dynamic;  Mori-Zwanzig; Krylov subspace; {\red Coarse grained molecular models}.}%Use showkeys class option if keyword
\maketitle

\section{\label{sec:intro}Introduction}

Direct numerical approaches  based on molecular interactions have become  standard
computational, as well as modeling,  tools nowadays for modeling molecular structures.
For dynamics problems, the trajectory of each atom can be described by the Newton's equations of motion,
\begin{equation}\label{eq: md}
\left\{
\begin{aligned}
 \dot{\bm x}_i &= \bm v_i,\\
m_i \dot{\bm v}_i &= -\frac{\partial V}{\partial \bm x_i} = \bm f_i(\bm x).
\end{aligned}\right.
\end{equation}
This approach is the essences of  the molecular dynamics (MD) modeling. The interatomic potential $V=V( x_1,  x_2, \cdots, x_{3N})$ embodies the interactions
between particles (atoms) through the changes of bond lengths, bond angles, dihedral angles, electrostatics, van der Waals
etc \cite{Schlick2002}.

Direct MD simulations  capture all the physics in a biological system, but they  particularly suited for studying small scale transitions due to the computational complexity.
Meanwhile, most biological processes are intrinsically multiscale: The overall dynamics consists of large number of atoms associated with
many different types of motions, spanning a wide range of time scales \cite{Schlick2002}.  In fact, typical biological functions begin at the $10^{-5}$s  time scale, which is far beyond the reach of direct MD simulations.

To overcome this significant modeling difficulty, much effort has been devoted to developing coarse-grained (CG) molecular models  to access processes on a longer time scale. Problems of this type have been identified as one of the most important and challenging problems in molecule modeling \cite{Leach01}. One of the key components in a CG model
is to find out the  direct interaction of the CG variables,  represented, \eg, by the many-body potential of mean force (PMF) \cite{baaden_coarse-grain_2013,gohlke_natural_2006,noid_perspective:_2013,noid_multiscale_2008,praprotnik_multiscale_2008,riniker_developing_2012,rudzinski_role_2012}.
In a CG approach, this interaction, in terms of forces,  can in principle be obtained by integrating out the remaining degrees of freedom  \cite{noid_perspective:_2013}.
However in practice, approximation schemes have to be introduced, and the main issue for PMF is to ensure the consistency   with the original full molecular interaction
as well as to control the accuracy.  We refer to the reviews \cite{noid_perspective:_2013,riniker_developing_2012} for the recent progress and existing issues.

The calculation of PMF is often formulated based on a thermodynamic consideration. In particular, one considers a system where the remaining degrees of freedom are at a conditional equilibrium. Another remarkable approach is through the generalized Langevin equations (GLE) , which can be derived directly from the equations of motion \eqref{eq: md} using the Mori-Zwanzig (MZ) projection formalism \cite{ChKaKu98,ChHaKu02,Mori65,Zwanzig73}.  The mode has been considered by many researchers over the years \cite{berkowitz1983generalized,bu2003vibrational,IzVo06,kamberaj2011theoretical,lange2006collective,oliva2000generalized,sagnella2000time,stepanova_dynamics_2007,Zwanzig61,Li2009c}. The MZ projection procedure, when the conditional expectation is used as projector, yields an averaged force, which is consistent with that in the PMF approach \cite{ChHaKu02,ChSt05,lange2006collective}. In addition, the formalism gives rise to a {\it history-dependent term}, which with reasonable approximations, simplifies to a linear convolutional term with a memory function, and a {\it random noise} term, which is consistent with the memory function via the
second fluctuation-dissipation theorem (FDT) \cite{Kubo66}.

The main practical difficulty in implementing the GLE is the computation of the memory function.  In some cases,
Markovian approximations can be made \cite{hijon2006markovian,kauzlaric2011bottom,kauzlaric2012markovian} to reduce the GLE to Langevin equation,
or one may simply use exponential functions \cite{oliva2000generalized}, assuming a rapid decay. However,
it is difficult to quantify and control the modeling error in such an ad hoc approximation.  A more systematic approach is to related the memory function to correlation functions, \eg, the velocity auto-correlation function (VACF), which is computed
from equilibrium MD simulations. For instance, Berkowitz et al  \cite{berkowitz1981memory} considered a GLE where
the mean force term is linear, and then derived an integral equation of Volterra type for the memory function. As input of the integral equation, the correlation function of the velocity and position are obtained from MD experiments.
This has been the approach followed by many other groups \cite{kamberaj2011theoretical,lange2006collective,sagnella2000time}.
 In general, the calculation of VACF tends to be expensive due to the large size of the system. More importantly, the sampling of the random noise is still a challenge.  In this paper, we propose a more efficient approach to obtain the memory functions without performing direct MD simulations. The method for computing the kernel functions is based on the Krylov-subspace method, motivated by the numerical methods for evaluating matrix functions. We will present the algorithm, and detailed implementation procedure. As will be shown, this approach offers the added advantage that the random force term can be approximated in the same subspace, and it automatically satisfy the second FDT.
 It is important to point out that the memory functions will depend on how the CG variables are selected, and what reduction procedure is used. The point will be illustrated and clarified using two reduction methods, and three different selection schemes for the CG variables.

The rest of the paper is organized as follows: We first discuss the reduction method of Mori-Zwanzig, from which we derive the exact expression of the memory functions. Then, we present an efficient numerical algorithm to compute these functions. Examples are given in the following section to demonstrate the effectiveness of the methods.

\section{\label{sec:GLE}The derivation of Generalized Langevin models}

The  generalized Langevin (GLE) models can be derived from many different coarse graining procedures, \eg,
by using appropriate linearization procedure \cite{stepanova_dynamics_2007}.
A more systematic procedure is the Mori-Zwanzig  projection formalism \cite{Zwanzig73,Mori65}. Here we will consider two different projection operators, and derive two types of GLEs models.   In particular, we derive an explicit expression for the memory function.

We start with the full molecular dynamics (MD) model,
\begin{equation}\label{eq: md0}
 m_i \ddot{\bm x}_i = \bm f_i(\bm x).
\end{equation}
Here $\bm x=(x_1, x_2, \cdots, x_{3N})$ denotes the position of {\it all} the atoms. Further, we let
$\bm v = \dot{\bm x}$ be the velocity.

Let us introduce a scaling,
\begin{equation}
 \bm x \to {m}^\half \bm x,\quad \bm v \to {m}^\half \bm v, \quad \bm f \to m^{-\half} \bm f.
\end{equation}
This reduces the equation \eqref{eq: md0} to
\begin{equation}\label{eq: md1}
  \ddot{\bm x} = \bm f(\bm x),
\end{equation}
which is expressed in a vector form. The coarse-graining procedure will be applied to these rescaled equations. In particular, the position will be mass weighted.

\smallskip

The collective motions are often represented in terms of the dynamics of a number of coarse-grained variables. We will define such variables through a projection to a subspace. Toward this end, we let $X=\mathbb{R}^{3N}$ be the entire configuration space, and $Y\subset X$ be a subspace with dimension $M$; $M\ll3N$.
Specific examples of such subspaces will be discussed later. To derive explicit formulas, let us choose
a set of orthonormal basis vectors of $Y$,  denoted here by  $(\varphi_1, \varphi_2, \cdots, \varphi_M).$ By grouping these vectors, we form a \(3N \times M\) matrix $\Phi$. Further, we let $(\psi_1, \psi_2, \cdots, \psi_{3N-M})$ be an orthonormal basis for the orthogonal complement of the subspace $Y$, denoted by $Y^\perp$. They form a \(3N \times (3N-M)\) matrix $\Psi$. In practice, it is often difficult, if not impossible, to construct the matrix $\Psi$.
Nevertheless, we will use this set of basis to express certain functions, and then we will discuss how to approximate these functions without actually computing $\Psi$.

To proceed, we define the CG variables through the projection to the subspace $Y$:
{\red \begin{equation}
\begin{aligned}
 \bm q= & \Phi^T \bm u,\\
 \bm p= & \Phi^T \bm v,
\end{aligned}
\end{equation}}
{\red where $\bm u= \bm x- \bm x^0$ is the displacement to the equilibrium state $\bm x^0$.
The displacement is often easier to work with, and we further switch the notation $\bm f(\bm x)$ to $\bm f(\bm u).$ }
Since all the columns  in $\Phi$ are unit vectors, $\bm q, \bm p \in \mathbb{R}^M$ can be regarded as average position and average velocity, respectively.  Similarly, one can define $\bm \xi= \Psi^T \bm u$ and $\bm \eta= \Psi^T \bm v$; $\bm \xi,\bm \eta \in \mathbb{R}^{3N-M}$. They represent the additional degrees of freedom, referred to as {\it under-resolved variables}, and they will not appear explicitly in the CG models.

It is clear now that for any $\bm u$ or $\bm v$, we have a unique decomposition in the form of,
{\red \begin{equation}
\begin{aligned}
  \bm u= &\Phi \bm q + \Psi \bm \xi,\\
  \bm v= &\Phi \bm p + \Psi \bm \eta.
\end{aligned}
\end{equation}}

The first step of the MZ reduction procedure is to express the time evolution of the CG variables. This is best represented by a semi-group operator, {\it i.e.,} for any dynamical variable $\bm y(t)$, we have $\bm y(t)= e^{t\cal{L}}\bm y$,  where
   the operator $\cal{L}$ is given by,
\begin{equation}
   \mathcal{L}= {\bm v}\cdot \frac{\partial}{\partial \bm x} + \bm f \cdot \frac{\partial}{\partial \bm v}.
   \end{equation}
As is customary in statistical mechanics theory, we use $\bm y$ to denote the initial value,  {\it i.e.,}
 $\bm y= \bm y(0)$, and these differential operators are defined with respect to the initial coordinate and momentum \cite{Balescu76,EvMo08,Zwanzigbook}. More specifically, the solution ($\bm x(t)$ and $\bm v(t)$) of the MD model \eqref{eq: md} at time $t$ depends on  the initial condition $\bm x$ and $\bm v$. Such dependence defines a symplectic mapping \cite{arnol1989mathematical}.   As a result, any dynamic variable $\bm y$, as a function of $\bm x(t)$ and $\bm v(t)$, are also functions of $\bm x$ and $\bm v$. The partial derivatives in $\cal{L}$ should be calculated with respect to the initial
condition.

In order to distinguish thermodynamic forces of different nature, one defines a projection operator $\cal{P}$, with its complementary operator given by $\cal{Q}=I-\cal{P}$.
It can either be defined as a projection to a subspace \cite{Mori65} or a conditional average \cite{Zwanzig73,ChKaKu98}. This will be discussed separately in the next section.

 Once the dynamic variables and the projection are defined, the Mori-Zwanzig procedure yields the effective model \cite{Mori65,Zwanzig73},
\begin{equation}\label{eq: MZ}
    \frac{d}{dt} \bm y(t) = e^{t{\cal L}} {\cal{PL}} \bm y + \int_0^t e^{(t-s){\cal L}} K(s) ds + R(t),
\end{equation}
where,
\begin{equation}\label{eq: R}
    R(t)= e^{t{\cal{QL}}} {\cal{QL}} \bm y,
\end{equation}
and,
\begin{equation}\label{eq: K}
 \quad K(t)={\cal{PL}}R(t).
\end{equation}

The first term on the right hand side of \eqref{eq: MZ} is typically considered as the reversible thermodynamic force. The second term represents the history dependence and provides a more general form of frictional forces. It dictates the strong coupling with the under-resolved variables. The last term, $R(t)$, takes into account the influence of the under-resolved variables, in the form of a random force.  Next, we discuss the specific forms of the memory function and the random noise for different choices of the projection operator.

\subsection{Orthogonal Projection}

Here we choose the following projection: For any function {\red $g(\bm u)$, or $g(\bm v),$} we define,
{\red \begin{equation}
 \cal{P} g(\bm u) = g(\Phi \bm q), \quad  \cal{P} g(\bm v) = g(\Phi \bm p).
\end{equation}}
The operator is a projection since $\Phi^T \Phi=I.$
This is motivated by the Galerkin method for coarse-graining MD models \cite{Li14}.

If $\bm y= \bm q,$ the MZ equation is reduced to,
\begin{equation}
 \frac{d}{dt} \bm q(t) = \bm p(t).
\end{equation}
No memory term arises from this equation.

Next, we let $\bm y=\bm p$. We will derive the CG model in several steps. First we start with the random noise $R(t)$.  At $t=0$,
we find $R= {\cal{QL}} \bm p(0)$ from \eqref{eq: R}.
\[ R=\cal{Q} \Phi^T \bm f = \Phi^T \Big[\bm f(\bm u)  - \bm f(\Phi \bm q) \Big].\]
In order to simplify this term, we introduce the approximation,
{\red \begin{equation}\label{eq: lin-approx}
\bm f(\bm u)  \approx \bm f(\Phi \bm q) - A \Psi \bm \xi.
\end{equation}}
{\blue In principle, one can choose $A=-\nabla \bm f( \Phi \bm q)$.}
{\blue But here we let } {\red $A=-\nabla f(\bm 0)$}, i.e., the hessian matrix of the potential energy at a local minimum $\bm x^0$, {\blue which has the same  second order accuracy of approximation near the reference position.}

With this approximation, we find that,
\( R \approx -\Phi^TA \Psi \bm \xi.\)
Applying the operator  $\cal{QL}$, we get,
\(\cal{QL} R \approx- \Phi^T A \Psi \bm \eta.\) We proceed to compute $(\cal{QL})^2 R$. A direct calculation yields,
{\red \(  (\cal{QL})^2 R \approx - \cal{Q} \Phi^T A \Psi \Psi^T \bm f(\bm u),\)}
which by a similar approximation \eqref{eq: lin-approx}, can be written as,
\(
    (\cal{QL})^2 R \approx \Phi^T A \Psi \Psi^T A \Psi \bm\xi.\) Similarly,
\[ (\cal{QL})^3 R \approx \Phi^T A \Psi \Psi^T A \Psi \bm\eta, \]
\[(\cal{QL})^4 R \approx -\Phi^T A \Psi \Psi^T A \Psi \Psi^T A \Psi \bm\xi.\]

Repeating such calculations, we find that the random noise can be approximated by
\begin{equation}
  R(t) \approx -\Phi^T A \Psi \Big[\cos\big(\Omega t\big)\bm \xi + \Omega^{-1} \sin \big(\Omega t\big) \bm \eta\Big],
\end{equation}
where $\Omega= \wh{A}^\half$ with
$\wh{A}=\Psi^T A \Psi$. This can be verified by examining the Taylor expansion of the trigonometric functions.

We now turn to the function $K(t)= \cal{PL} R(t).$ With the approximation of $R(t)$, we obtain,
\begin{equation}
  K(t) \approx -\Phi^T A \Psi \Omega^{-1} \sin \big(\Omega t\big) \Psi^T \bm f(\Phi \bm q).
\end{equation}
To further simplify this, we make another approximation that $\bm f(\Phi \bm q) \approx - A \Phi \bm q$ in this expression
\footnote{This linearization is again the reference position, in which case the matrix is
the hessian of the potential energy.
This seems to be the only way to obtain a linear convolutional form of the memory function. }
 , which leads to,
\[ K(t) \approx \Phi^T A \Psi \Omega^{-1} \sin \big(\Omega t\big) \Psi^T A \Phi \bm q. \]
This simplifies the integral to a convolutional form,
\[ \int_0^t \beta(s) \bm q(t-s) ds,\]
where the $M\times M$ matrix function $\beta$ is given by,
\begin{equation}
\beta(t) =  \Phi^T A \Psi \Omega^{-1} \sin \big(\Omega t\big) \Psi^T A \Phi.
\end{equation}

Collecting terms, we obtain the GLE,
\begin{equation}\label{eq: gle1}
 \ddot{\bm q}= F(\bm q) + \int_0^t \beta(t-s) \bm q(s)ds + R(t).
\end{equation}

The first term in the GLE \eqref{eq: gle1} is related to the inter-molecular force as follows:
\begin{equation}
F(\bm q)=\Phi^T \bm f(\Phi \bm q).
\end{equation}

\subsection{Projection via Conditional Expectation}
Another choice of the projection operator is the conditional expectative, which for the canonical ensemble,
is given by,
\begin{equation}
\begin{split}
&\cal{P} g(\bm x, \bm v)= E\Big[ g(\bm x, \bm v)| \bm q, \bm p\Big]\\
&\stackrel{\text{def}}{=}\frac{\dsp\int_{\mathbb{R}^{6N}} g(\bm x,\bm v)e^{-\beta[V(\bm x)+\half\bm v^2]}\delta(\bm q-\Phi^T\bm x)\delta(\bm p-\Phi^T \bm v)d\bm xd\bm v}{\dsp\int_{\mathbb{R}^{6N}}e^{-\beta[V(\bm x)+\half \bm v^2]}\delta(\bm q-\Phi^T \bm x)\delta(\bm p-\Phi^T\bm v) d\bm xd\bm v}
\end{split}
\end{equation}
Here $\beta=\frac{1}{k_B T}$ is the inverse temperature, and the delta functions are introduced to enforce the given conditions.

Again we start with the construction of the random noise in the MZ equation \eqref{eq: MZ}. Here we introduce two approximations. First, we let $A\approx -\nabla \bm f$ be an approximate hessian of the potential energy, and we approximate the projection by,
\begin{equation}\label{eq: approx1}
\begin{split}
&\cal{P} g(\bm x, \bm v)\approx \\
&\frac{\dsp\int_{\mathbb{R}^{6N}} g(\bm x, \bm v)e^{-\beta[\half \bm x^T A \bm x + \half \bm v^2]} \delta(\bm q - \Phi^T \bm x) \delta(\bm p - \Phi^T \bm v)   d\bm x d\bm v}
{\dsp\int_{\mathbb{R}^{6N}} e^{-\beta[\half \bm x^T A \bm x + \half \bm v^2]} \delta(\bm q - \Phi^T \bm x) \delta(\bm p - \Phi^T \bm v)   d\bm x d\bm v}.
\end{split}
\end{equation}
 As a result, the expectation is with respect to a multi-variant Gaussian distribution.

The second approximation also involves the same linearization used in the previous section,
\begin{equation}\label{eq: approx2}
 \bm f(\bm x) \approx \bm f(\Phi \bm q) - A \Psi \bm \xi.
\end{equation}
To facilitate the following calculations, we define projection matrices \cite{yanai2011projection,Li2009c},
\begin{equation}
 \begin{aligned}
 P_v= & \Phi \Phi^T,\quad  Q_v= I - P_v= \Psi \Psi^T, \\
 P_x= & A^{-1} \Phi \big(\Phi^T A^{-1} \Phi\big)^{-1} \Phi^T, \quad Q_x= I-P_x.
 \end{aligned}
\end{equation}
In particular, we have $\cal{P} \bm v= P_v\bm v$, and with the approximation \eqref{eq: approx1}, we have,
\[ \cal{P} \bm x\approx P_x \bm x.\]
Therefore, the projection operator has been turned into a matrix-vector multiplication.

The following identities can be easily verified,
\begin{align}\label{eq: ids}
& P_v P_x= P_v, \\ & Q_x Q_v = Q_v, \\ &Q_x P_v=Q_x - Q_v = P_v - P_x.
\end{align}

We proceed to compute the random noise. At $t=0$, $\cal{QL}\bm p= \cal{Q} \Phi^T \bm f(\bm x).$
By invoking the two approximations, we find that,
\begin{equation}
\begin{split}
R&=\cal{QL}\bm p\\&\approx \cal{Q} \Phi^T \big(\bm f(\Phi \bm q) - A Q_v \bm x\big) \\&\approx
-\Phi^T A Q_v Q_x \bm x\\&= -\Phi^T A \Psi \Psi^T Q_x \bm x.
\end{split}
\end{equation}
In addition, we have,
$$\cal{QL}R \approx -\Phi^T A \Psi \Psi^T \bm v.$$
Repeating these steps, we have,
\[
\begin{aligned}
(\cal{QL})^2R & \approx \Phi^T A \Psi \Psi^T A \Psi \Psi^T Q_x \bm x=\Phi^T A \Psi \widehat{A} \Psi^T Q_x \bm x,\\
(\cal{QL})^3 R & \approx \Phi^T A \Psi \widehat{A} \Psi^T \bm v,\\
(\cal{QL})^4 R & \approx \Phi^T A \Psi \widehat{A}^2 \Psi^T Q_x \bm x,
 ... &
 \end{aligned} \]
Again we defined $\wh{A}=\Psi^T A \Psi$. These calculations suggest that the random noise may be approximated by,
\begin{equation}\label{eq: R(t)}
\begin{split}
&R(t) \approx\\ &-\Phi^T A \Psi \Big[ \cos(\Omega t) \Psi^T Q_x \bm x + \Omega^{-1}\sin (\Omega t)
\Psi^T \bm v\Big],
\end{split}
\end{equation}
which can be validated by checking each term in the Taylor series.

With the approximation of $R(t)$, we can approximate $K(t)$ by,
\begin{equation}
\begin{aligned}
K(t)\approx& \Phi^T A \Psi \cos (\Omega t) \Psi^T A^{-1} \Phi (\Phi^T A^{-1} \Phi)^{-1}\bm p\\
&- \Phi^T A \Psi \Omega^{-1}\sin (\Omega t) \Psi^T \big(\bm f(\Phi \bm q) +  A \Phi \bm q\big).
\end{aligned}
\end{equation}
Here we have used the first and third identities in \eqref{eq: ids}.

Similar to the previous section, we neglect the second term using \eqref{eq: approx2} \footnote{ Otherwise
a nonlinear convolution term appear. Such an approximation is
appropriate when the system is around the reference position. But the error would be large during structural change, in which case it would be necessary to keep this term.}.
 As a result, we obtain a memory function,
\begin{equation}
 \theta(t) =-\Phi^T A \Psi \cos (\Omega t) \Psi^T A^{-1} \Phi (\Phi^T A^{-1} \Phi)^{-1}.
\end{equation}
Further, the memory term is reduced to a convolutional integral,
\begin{equation}
 -\int_0^t \theta(s) \bm p(t-s)ds.
\end{equation}
Notice that the memory function involves the coarse-grained momentum instead of the coarse-grained coordinate.

Using the matrix identity \cite{Li2009c},
\begin{equation}
 A^{-1} \Phi (\Phi^T A^{-1} \Phi)^{-1} = \Phi - \Psi (\Psi^T A \Psi)^{-1} \Psi^T A \Phi,
\end{equation}
we can simplify the memory function to,
\begin{equation}
 \theta(t) =\Phi^T A \Psi \cos (\Omega t) \Omega^{-2} \Psi^T A \Phi. \label{kernel}
\end{equation}

To get some insight, we let the eigenvalues of $\wh{A}$ be $\lambda_i$, and let $\bm w_i$ be the associated eigenvectors. Then, we can express the kernel function as follows,
\begin{equation}
 \theta(t)= \sum_{\lambda \in \sigma(\wh{A})} \frac{1}{\lambda_i}\cos(\sqrt{\lambda_i} t)
  (\Phi^T A\Psi_i \bm w_i)\otimes (\Phi^T A\Psi \bm w_i).
\end{equation}
Further, let $\bm r_i= A\Psi_i \bm w_i - \lambda_i \Psi_i \bm w_i$. This can be interpreted as
the residual error, when $\lambda_i$ is viewed as the approximate eigenvalue of $A$ obtained by a projection to the orthogonal complement. A direct substitution yields,
\begin{equation}\label{eq: theta}
 \theta(t)= \sum_{\lambda \in \sigma(\wh{A})} \frac{1}{\lambda_i}\cos(\sqrt{\lambda_i} t)
  (\Phi^T \bm r_i)\otimes (\Phi^T \bm r_i).
\end{equation}
Intuitively, when the eigenvalues are well approximated within the initial subspace $Y$, they make
less contribution to the memory function.

With the condition expectation chosen as the projection operator, the first term in the MZ equation \eqref{eq: MZ} has a
natural interpretation. To explain this, we define the free energy by integrating out the under-resolved variables,
\begin{eqnarray}
&W(\bm q, k_BT) = -k_B T \text{ln} Z\\&Z= \int_{\mathbb{R}^{3N}} e^{-\beta V(\bm x) } \delta(\bm q - \Phi^T \bm x)   d\bm x.
\end{eqnarray}
Then the first term in \eqref{eq: MZ} coincides with the mean force $-\nabla_{\bm q} W(\bm q, k_BT)$.

Now we can collect all the terms and the GLE is expressed as,
\begin{equation}\label{eq: gle2}
\ddot{\bm q}= -\nabla_{\bm q} W(\bm q, k_BT) - \int_0^t \theta(t-s)\dot{\bm q}(s)ds + R(t).
\end{equation}

With the approximation of the probability density, we see that \(\bm \xi\) and \(\bm \eta\)
follow the conditional distribution,
\begin{equation}
 (\bm \xi, \bm \eta) \sim e^{-\beta\Big[\bm \zeta^T\Psi^TA\Psi
 \bm \zeta+\bm\eta^T\Psi^T\Psi\bm\eta\Big]},
\end{equation}
where $\bm \zeta=\Big(\bm\xi+(\Psi^TA\Psi)^{-1}\Psi^TA\Phi^T\bm q\Big)$.

In addition, we have \[Q_x\bm x= \Psi \big( \bm \xi + (\Psi^T A \Psi)^{-1} \Psi^T A \Phi^T \bm q \big).\]
Therefore, the random process $R(t)$ in \eqref{eq: R(t)} is a Gaussian  process. Furthermore, with direct calculation, we can verify that it is stationary with zero mean and it
satisfies the second fluctuation-dissipation theorem (FDT) \cite{Kubo66},
\begin{equation}
 \Big\langle R(t) R(s)^T \Big\rangle= k_B T \theta(t-s).
\end{equation}
Based on the theory of Gaussian processes \cite{Doob44}, $R(t)$ is uniquely determined by the correlation function. Thus, the GLE is closed.  The FDT a critical property of the generalized Langevin model. It is a necessary condition to ensure
that the system will approach to a thermodynamic equilibrium   \cite{Kubo66}. Therefore, it is also important to preserve this
condition at the level of numerical approximations. This will be discussed in the next section.

In contrast, the random noise derived from the previous section is not stationary. However, notice that
\begin{equation}
\dot{\theta}(t)= -\beta(t).
\end{equation}
Using integration by parts, one can show that the memory functions and random noises in the GLEs \eqref{eq: gle1} and \eqref{eq: gle2} can be related to one another. For the rest of the paper, we will focus on the GLE \eqref{eq: gle2} and
the memory function $\theta(t)$. The function $\beta(t)$ can be computed using a similar procedure.

\section{\label{sec:krylov} A Krylov subspace approximation of the kernel function}
In most of previous works, the memory functions are computed from molecular dynamics simulations.
%In particular, by neglecting the average force, the generalized Langevin equation can be reduced to,
%\begin{equation}
% \ddot{\bm q} = - \int_0^t \theta(t-s) \dot{\bm q}(s)ds + R(t).
%\end{equation}
%Assume that $\langle R(t) \bm p(0)^T \rangle=0.$ By multiplying the equation by $\bm p(0)^T$, and assuming that $\bm p(t)$ is a stationary process, one can derive an integral equation,
%\begin{equation}
% \frac{d}{dt} C(t) = - \int_0^t \theta(t-s) C(s) ds.
%\end{equation}
%where $C(t)$ is the velocity correlation function.
%\begin{equation}
%C(t)= \langle \bm p(t) \bm p(0)^T \rangle.
%\end{equation}
%For systems of moderate size, $C(t)$ can be sampled from a full MD simulation. Then, one may solve the integral equation to obtain the memory function.
%
%
%However, it is clear that the memory function will depend on the force term. To demonstrate such dependence, we consider the follow two scalar
%GLEs.
%\begin{equation}
%  \ddot{q}= - \int_0^t \theta(t-s)\dot{q}(s) ds + R(t),
%\end{equation}
%and
%\begin{equation}
%  \ddot{q}= -k q - \int_0^t \theta(t-s)\dot{q}(s) ds + R(t).
%\end{equation}
%If we let $C(t)$ be the velocity auto-correlation and $\wt{C}(\lambda)$ be the corresponding Laplace transform. Then, for the first model, we have
%\begin{equation}
% \wt{\theta}(\lambda) = \frac{ k_BT - \lambda \wt{C}(\lambda)}{\wt{C}(\lambda)},
%\end{equation}
%while for the second model, we find that,
%\begin{equation}
%  \wt{\theta}(\lambda) = \frac{\lambda k_BT - (k/\lambda+\lambda) \wt{C}(\lambda)}{\wt{C}(\lambda)}.
%\end{equation}
%Clearly the kernel function depends on $k$.
%
%\smallskip
In this paper, we present another approach, based on the analytical expression of the kernel \eqref{eq: theta}.
Due to the matrix function form, we will use the Krylov subspace approximation, a popular method for computing
matrix functions \cite{saad_analysis_1992,diele_error_2009}. Next, we explain the general idea, and address some
implementation issues.

\subsection{Approximation using the Krylov spaces}

We first consider the approximation of $\theta(0)$ to illustrate the idea. Recall that $\wh{A}=\Psi^T A \Psi$,
and so $$\theta(0)=\Phi^T A \Psi \wh{A}^{-1} \Psi^T A \Phi.$$ Consider the vector $\bm b= \Psi^T A \varphi_k$ for some $k$, $1\le k\le M$, and  we define the Krylov subspace with {\it order} $m$,
\begin{equation}\label{eq: krylov1}
  K_m\big(\wh{A},\bm b\big)= \text{span}\Big\{\bm b, \wh{A} \bm b, \cdots, \wh{A\,}^m \bm b \Big\}.
\end{equation}
With the standard Lanczos algorithm \cite{Saad}, we can construct orthogonal basis vectors $B_m=[\bm b_1, \bm b_2, \cdots, \bm b_m]$ for $K_m\big(\wh{A},\bm b\big)$.  Further, it reduces the matrix $\wh{A}$ to the form,
\begin{equation}
\wh{A} B_m= B_m T_m + \gamma_{m+1} \bm b_{m+1} \bm \epsilon_m^T.
\end{equation}
The last term, which is a rank-one matrix, contains the error.

As a result, we make the approximation
\begin{equation}
 \wh{A}^{-1} \bm b \approx ||\bm b||_2 B_m T_m^{-1} \bm e_1.
\end{equation}
Therefore, the $(k,k)$ entry of $\theta(0)$ can be approximated by,
\begin{equation}
\theta_{k,k}(0)= \bm b^T\wh{A}^{-1} \bm b \approx ||\bm b||^2_2\bm e_1^T T_m^{-1} \bm e_1.
\end{equation}
The vector $\bm e_1$ is the standard basis vector. Consequently, the computation of the inverse of a large matrix is reduced to the inversion of a much smaller, tri-diagonal, matrix $T_m$ \cite{Saad}.

For the present problem, several issues arise:
\begin{enumerate}
\item {\red Both $||\bm b||^2_2$ and the matrix $\wh{A}$ are difficult to compute directly}, since the basis functions $\psi_i$  are usually not available;
\item There are a number of basis vectors $\varphi_i$ to begin with, {\blue and we need to compute the entire
matrix $\theta(t)$. The standard Krylov space method has to be implemented multiple times to obtain the entire matrix.}
\end{enumerate}

To overcome the first difficulty, we introduce a mathematically equivalent procedure based on the following observation.
Recall that  $Q_v=\Psi\Psi^T$, and we now define $\wt{A}= Q_vAQ_v$, and $\bm w=  A \varphi_k$. It can be
directly verified that,
\begin{equation}
\Psi K_m \big(\wh{A},\bm b\big)= \text{span}\Big\{\bm w, \wt{A}\bm w, \cdots, \wt{A}^m\bm  w\Big\}= K_m   \big(\wt{A},\bm w\big).
\end{equation}

In addition,
\begin{rem}
The Lanczos algorithm, when applied to the subspace $\Psi K_m \big(\wh{A},b\big)$, yields the same results
as those obtained from the  Lanczos algorithm applied to the subspace $K_m \big(\wt{A},w\big).$
\end{rem}

Now in the Krylov space $K_m \big(\wt{A},w\big)$, $\Psi$ is not involved. Further, $Q_v=I-P_v = I-\Phi\Phi^T$. This can be
drastically simplified when the basis functions in $\Phi$ are localized. One such example is the rotational-translational block
method (RTB) \cite{licu02,TaGaMaSa00},  which divides the entire molecule into non-overlapping blocks. In each block, the  rotational and translational degrees of freedom can be selected as basis functions. The explicit formulas can be found in
\cite{DeMi12,Wilson80}. When such basis functions are used, the matrix $P_v$ is block diagonal, and the matrix $Q$ can be easily computed. In fact,  the implementation of the above algorithm only involves
the product of $Q$ with another vector. The multiplication can be done separately in each block.

%Finally, we may make the approximation,
%\begin{equation}
% \Theta_{jk}(0)= \varphi_j^t A W_m T_m^{-1} e_1.
%\end{equation}
%and,
%\begin{equation}
% \Theta_{jk}(t)= \varphi_j^t A W_m T_m^{-1} \cos (T_m^\half t) e_1.
%\end{equation}

\smallskip

To address the second issue, we employ the {\it block} Krylov method and {\it block} Lanszos method. The application of
the block Krylov method can be found in \cite{ye_adaptive_1996}. Here we provide some details.

We  first let $V=A \Phi,$ and define,
\begin{equation}\label{eq: krylov2}
K_m\big(\wh{A},V\big)= \text{span}\Big\{V, \wh{A} V, \cdots, \wh{A}^m V \Big\}.
\end{equation}
The right hand side is interpreted as the linear combination of the columns of the matrices. It is a natural generalization of
the Krylov space \eqref{eq: krylov1}. To obtain orthonormal basis for the subspace, we follow the steps below:\\
\noindent{\bf Algorithm.} (Block Lanczos) Set $V_0=0$, $Z_0=V$ and $p_0=M$. For $j=1, 2, \cdots, m$, repeat:
\begin{itemize}
\item []{\bf Step 1.} Rank revealing QR factorization of the $n\times p_{j-1}$ matrix $Z_{j-1}$: $Z_{j-1} = Q_j R_{j-1}$. $R_{j-1}$ may be a permuted upper triangular matrix.
\item []{\bf Step 2.} Let $p_j=\text{rank}(Z_{j-1}),$ $V_j$ be the first $p_j$ columns of $Q_j$,
and $B_{j-1}$ be the first $p_j$ rows of $R_{j-1}$;
\item []{\bf Step 3.} $Z_j \longleftarrow A V_j - V_{j-1}B_{j-1}^T$;
\item []{\bf Step 4.} $A_j \longleftarrow V_j^T Z_j$;
\item []{\bf Step 5.} $Z_j \longleftarrow Z_j - V_j A_j$.
\end{itemize}

Let $E_1=[I \; 0 \; \cdots\; 0]^T$. We then have
\begin{equation}
\theta(0) \approx B^T_0 E_1^T T_m^{-1} E_1 B_0, B_0= (V^TV)^\half.
\end{equation}
Similarly, we have,
\begin{equation}
 \theta(t)\approx \wh{\theta}(t) \stackrel{def}{=} B^T_0 E_1^T T_m^{-1} \cos (T_m^\half t) E_1 B_0.\label{approx3}
\end{equation}
\subsection{Approximation of the random noise}

We now turn to the random noise $R(t)$, which can also be sampled within the Krylov subspace. More precisely, we state that,
\begin{rem}
Let $\wh{R}(t)$ be given by,
\begin{equation}
\wh{R}(t)= E_1^T \cos (T_m^\half t) \wh{\bm \xi} + E_1^T \sin(T_m^\half t)T^{-\half} \wh{\bm \eta},
\end{equation}
where $\wh{\bm \xi}$ and $\wh{\bm \eta} $ are independent normal random variables with zero mean and
variance $k_B T T_m^{-1}$ and $k_B T I$, respectively, then $\wh{R}(t)$ is stationary random noise with zero mean
and the correlation is given by,
\begin{equation}
 \Big\langle \wh{R}(t) \wh{R}(s)^T \Big\rangle=k_B T \wh{\theta}(t-s).
\end{equation}
\end{rem}
As a result, the sampling of the random force is reduced to the sampling of low-dimensional quantities $\wh{\bm \xi}$ and $\wh{\bm \eta} $. More importantly, the approximate random force $\wh{R}(t)$ and memory function $\wh{\theta}(t)$
still satisfy the fluctuation-dissipation theorem.

%Further, with the projection to the Krylov subspace, we can embed the GLE into a Markovian system. In particular, we have
%\begin{rem}
% The GLE model is equivalent to,
% \begin{equation}
%\left\{
%\begin{aligned}
% \ddot{\bm q} &= F(\bm q) - B_0^T E_1^T \bm z, \\
% \ddot{\bm z} &= -E_1 B_0  \bm q - T_m\bm z.
%\end{aligned}
%\right.
%\end{equation}
%\end{rem}
%This can be proved by solving the second equation, and then substitute it into the first equation, which will recover the GLE.

\section{\label{sec:eg} Examples}
\begin{figure}[tph]
\begin{center}
\caption{Cartoon picture of the structure of protein HIV-1 protease (PDB id:1DIF).}\label{pdbfig}
\end{center}
\end{figure}
In this section, we present some numerical results.
As an example, we choose a HIV-1 protease whose PDB id is 1DIF.
The protein contains 198 residues and 3128 atoms. The cartoon picture of the structure is shown in Fig. \ref{pdbfig}.

The kernel functions depend on the choice of the coarse-grained variables. In particular, it depends on the initial subspace.
Here, three different subspaces are considered:
\begin{itemize}
  \item {\bf Subspace-I:} The subspace spanned by the RTB basis corresponding to the translations and rotations of rigid blocks. The partition of the blocks is obtained from the partition scheme FIRST \cite{JaRaKu01}. The implementation was done by using the software PROFLEX.  The dimension of the subspace is 380.
  \item {\bf Subspace-II:}  The subspace generated by the RTB basis functions with each residue as a rigid block. There are 1188 basis functions in total.
  \item {\bf Subspace-III:} The subspace spanned by 540 low frequency modes, obtained from the principle component analysis (PCA) \cite{Tatsuoka}.  To obtain the basis functions,  trajectories are generated from direct molecular dynamics  simulations. These basis functions may not be localized. Nonetheless, we still choose this subspace due to its importance
  in dimension reduction.
\end{itemize}
For each subspace, we use the Krylov subspace methods and compute the approximate memory functions in (\ref{approx3}). For comparison, we also computed the exact memory function \eqref{kernel} using brutal force. {\blue The kernel functions have the unit of $\rm{eV/\AA^2}$.}

\begin{figure}[htbp]
\begin{center}
\includegraphics[height=5cm, width=9cm]{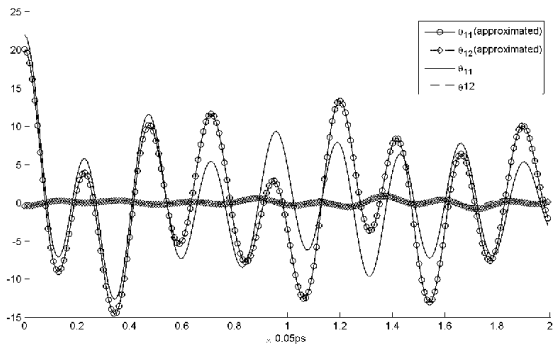}
\caption{Profiles of the kernels function for subspace-I: The first two entries $\theta_{11}(t)$ and $\theta_{12}(t)$
of the exact kernel function (lines without markers) produced by brutal-force computation according to (\ref{kernel}) and approximated kernel (\ref{approx3}) using the Krylov space method (lines with markers) with order 4. {\red These two entries are corresponding to the correlations  of the noises in the first two translational modes of the first rigid block.\label{fig1}}}
\end{center}
\end{figure}

\begin{figure}[htbp]
\begin{center}
\includegraphics[height=5cm, width=9cm]{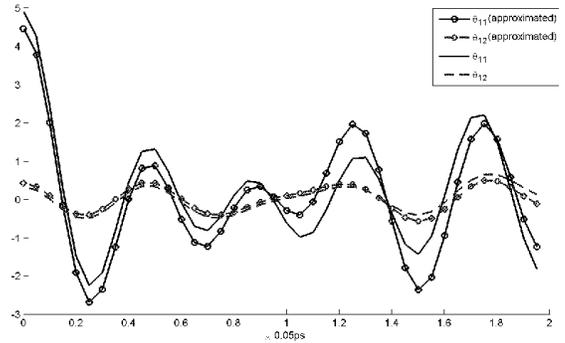}
\caption{Profiles of the kernels function for subspace-II:  The first two entries of the exact kernel function (lines without markers) produced by directly computation according to (\ref{kernel}) and the approximated kernel (\ref{approx3}) using the Krylov space method and subspace-II (lines with markers). {\red These two entries are corresponding to the correlations of the noises in the first two translational modes of the first rigid block.} The order of the Krylov space is 4.}\label{fig2}
\end{center}
\end{figure}

\begin{figure}[htbp]
\begin{center}
\includegraphics[height=5cm, width=9cm]{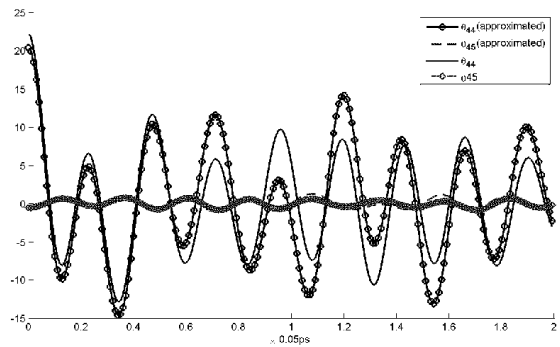}
\caption{Profiles of the kernels function for subspace-I: The first two entries of the exact kernel function (lines without markers) produced by brutal-force computation and the approximated kernel (\ref{approx3}) using Krylov space method (lines with markers) with order 4. {\red These two entries are corresponding to the correlations of the noises in the first two rotational modes of the first rigid block.}}\label{fig3}
\end{center}
\end{figure}

\begin{figure}[htbp]
\begin{center}
\includegraphics[height=5cm, width=9cm]{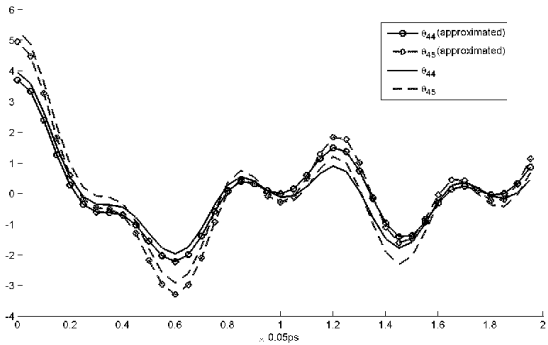}
\caption{Profiles of the kernels function for subspace-II: Two entries of the kernel function (lines without markers) produced by directly computation and the the approximated kernel using Krylov space method (lines with markers). {\red These two entries are corresponding to the correlations of the noises in the first two rotational modes of the first rigid block.} The order of the Krylov space is 4. }\label{fig4}
\end{center}
\end{figure}

\begin{figure}[htbp]
\begin{center}
\includegraphics[height=5cm, width=9cm]{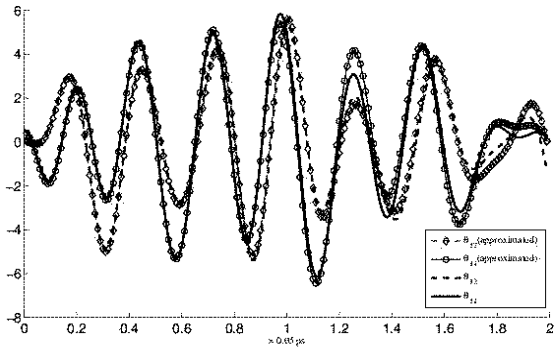}
\caption{Profiles of the kernels function for subspace-III: The first two entries of the exact kernel function produced by directly computation (lines without markers) and the approximated kernel using the Krylov space method (lines with markers) using subspace-III.  The order of the Krylov space  is 4.}\label{fig5}
\end{center}
\end{figure}

\begin{figure}[htbp]
\begin{center}
\includegraphics[height=5cm, width=9cm]{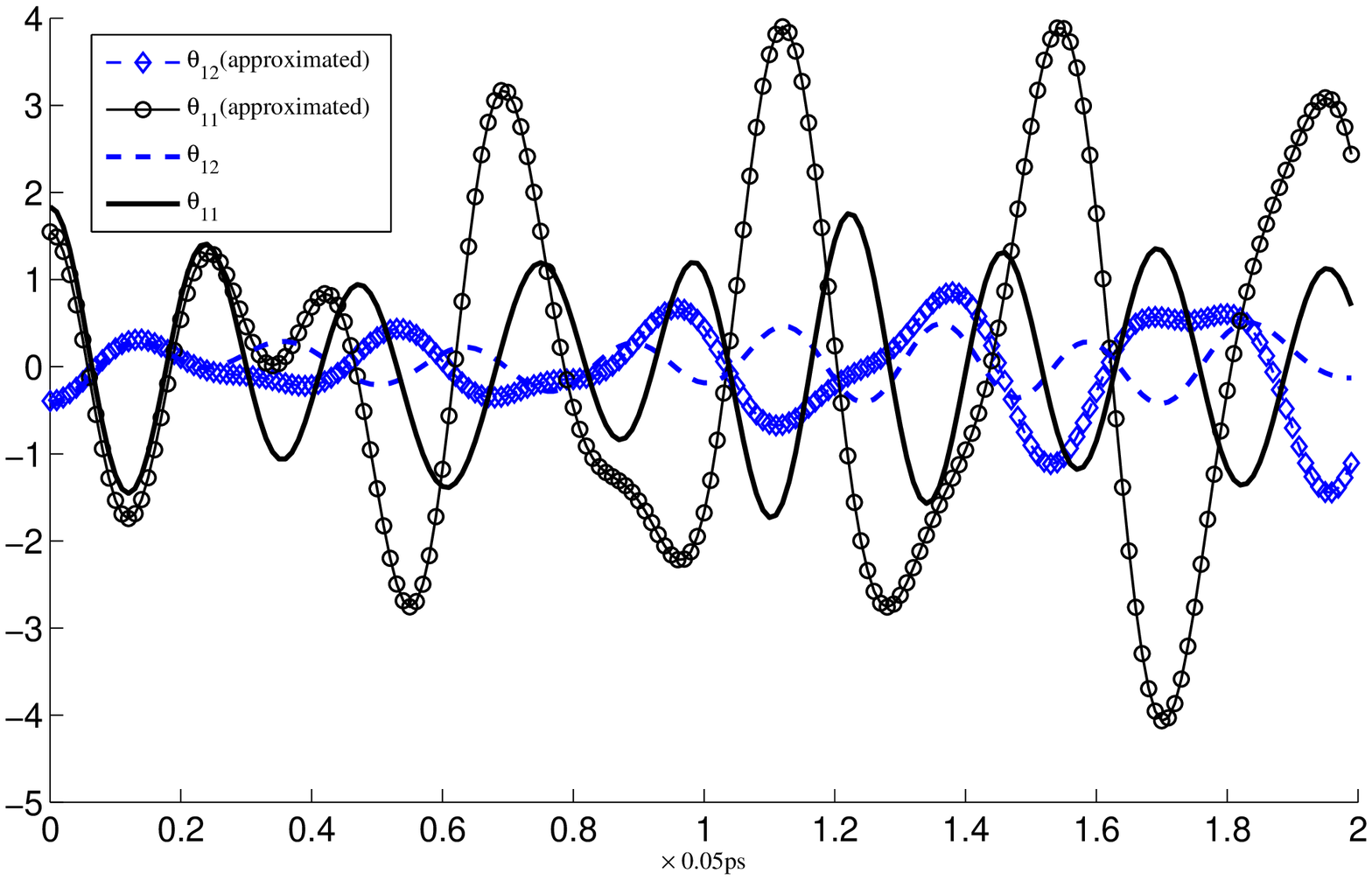}
\includegraphics[height=5cm, width=9cm]{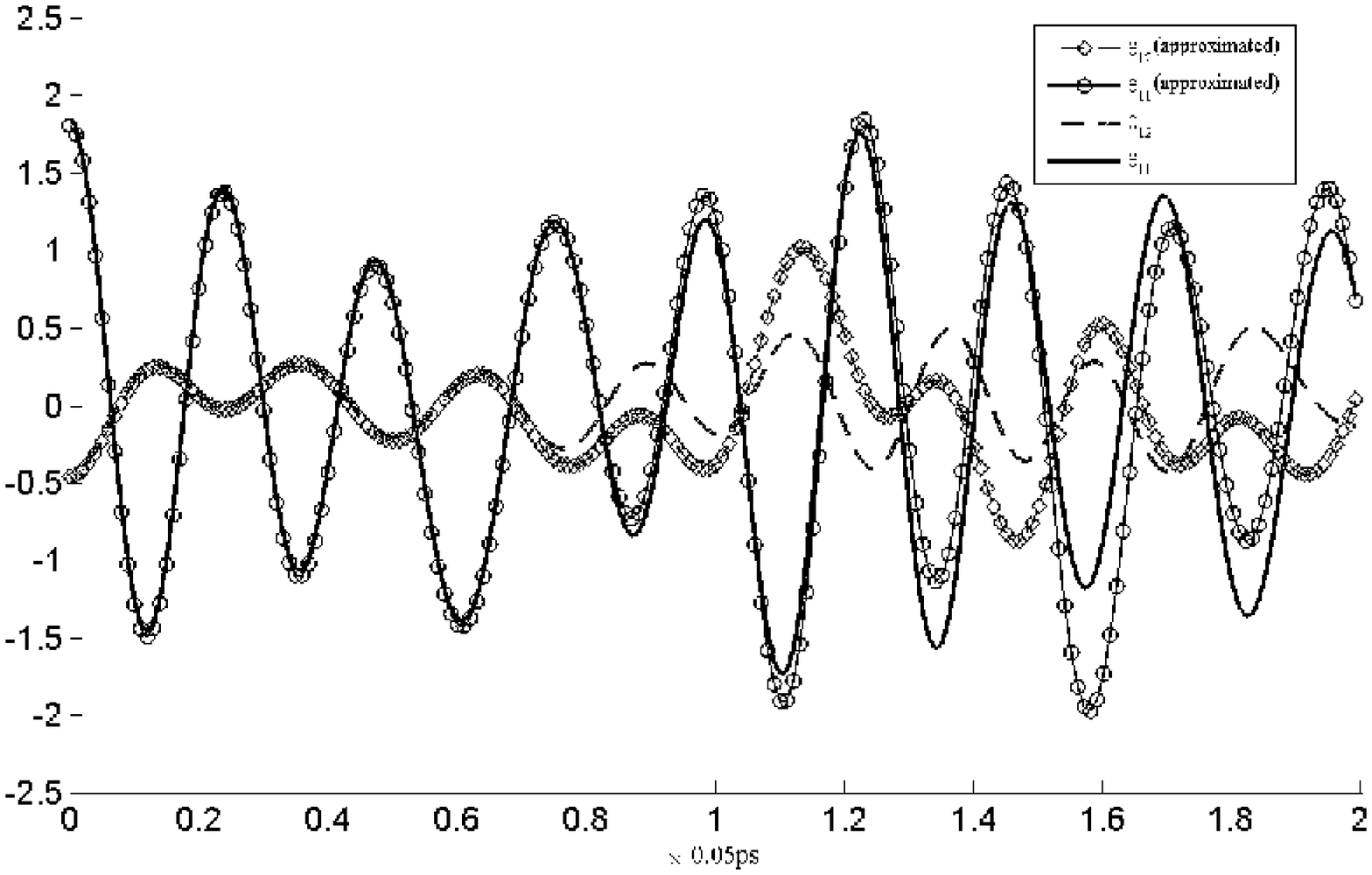}
\caption{Profiles of the first two entries of the exact kernel function produced by directly computation (lines without markers) and the approximated kernel using the Krylov method (lines with markers) for subspace-I.  Top:  Krylov space with order 2; Bottom: Krylov space with order 6.}\label{fig6}
\end{center}
\end{figure}

In Fig. \ref{fig1} - \ref{fig5}, {\red we show the profiles of the entries $\theta_{11}(t)$, $\theta_{12}(t)$
and $\theta_{44}(t)$, $\theta_{45}(t)$ } of the kernel function $\theta(t)$ within a time period of 0.1ps obtained from different computational methods and different coarse grained subspaces.
Based on these figures, we can see that the Krylov space method produces good approximations of the kernel functions,  especially at the beginning period.  Another observation is that these memory functions do not exhibit fast decay at this scale. Instead, they exhibit many oscillations, which indicate that a Markovian or exponential approximation is premature.
 Currently the order of the Krylov subspace in these examples are  4. If we increase the order of the Krylov subspace, the approximations will further improve, see Fig. \ref{fig6}.

 Fig. \ref{fig2}  also indicates that the memory functions for the residue-based subspaces look smoother. This is because the residue-based subspaces admit more low frequency modes than those of rigid bodies from the partitions of PROFLEX.

\begin{figure}[htbp]
\begin{center}
\includegraphics[height=5cm, width=9cm]{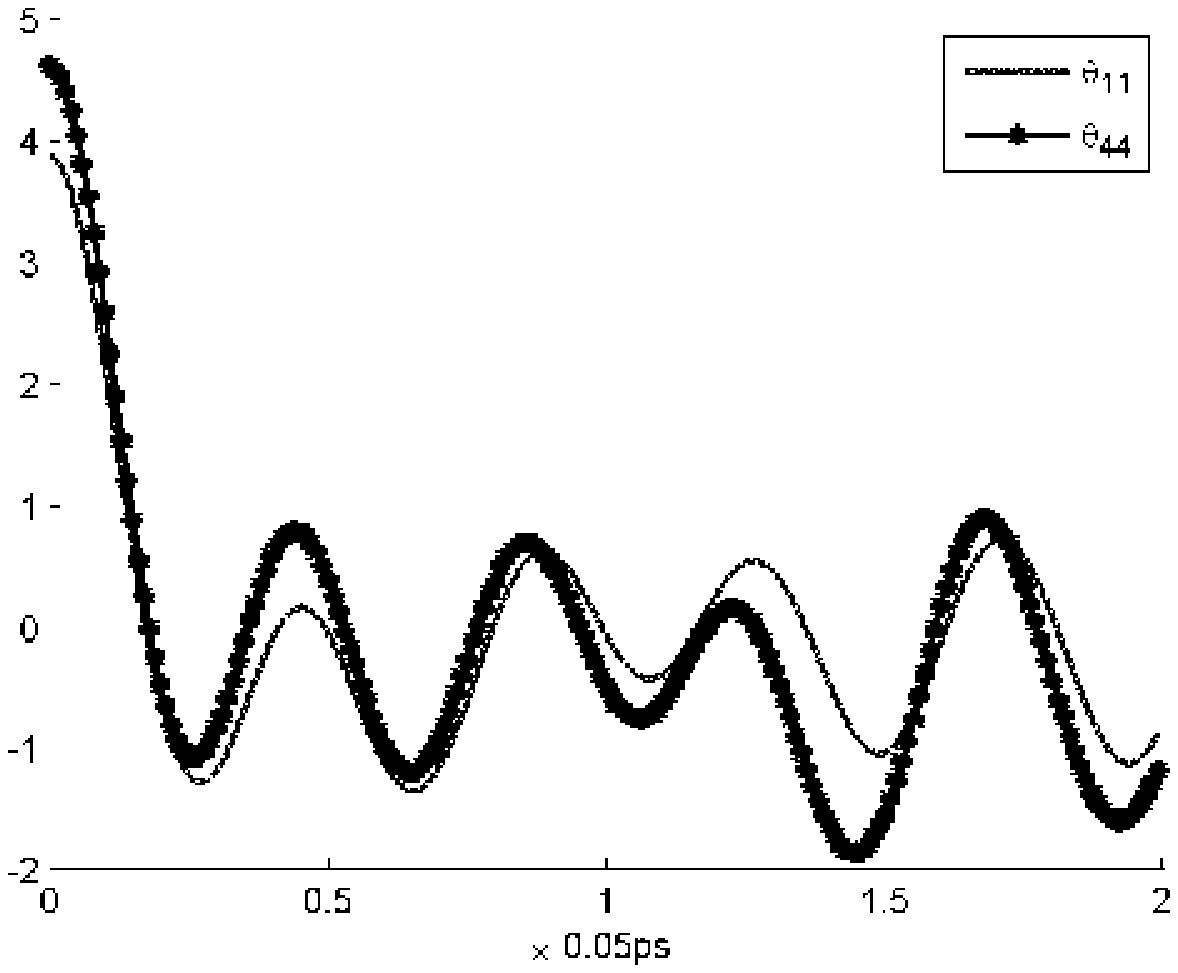}
\includegraphics[height=5cm, width=9cm]{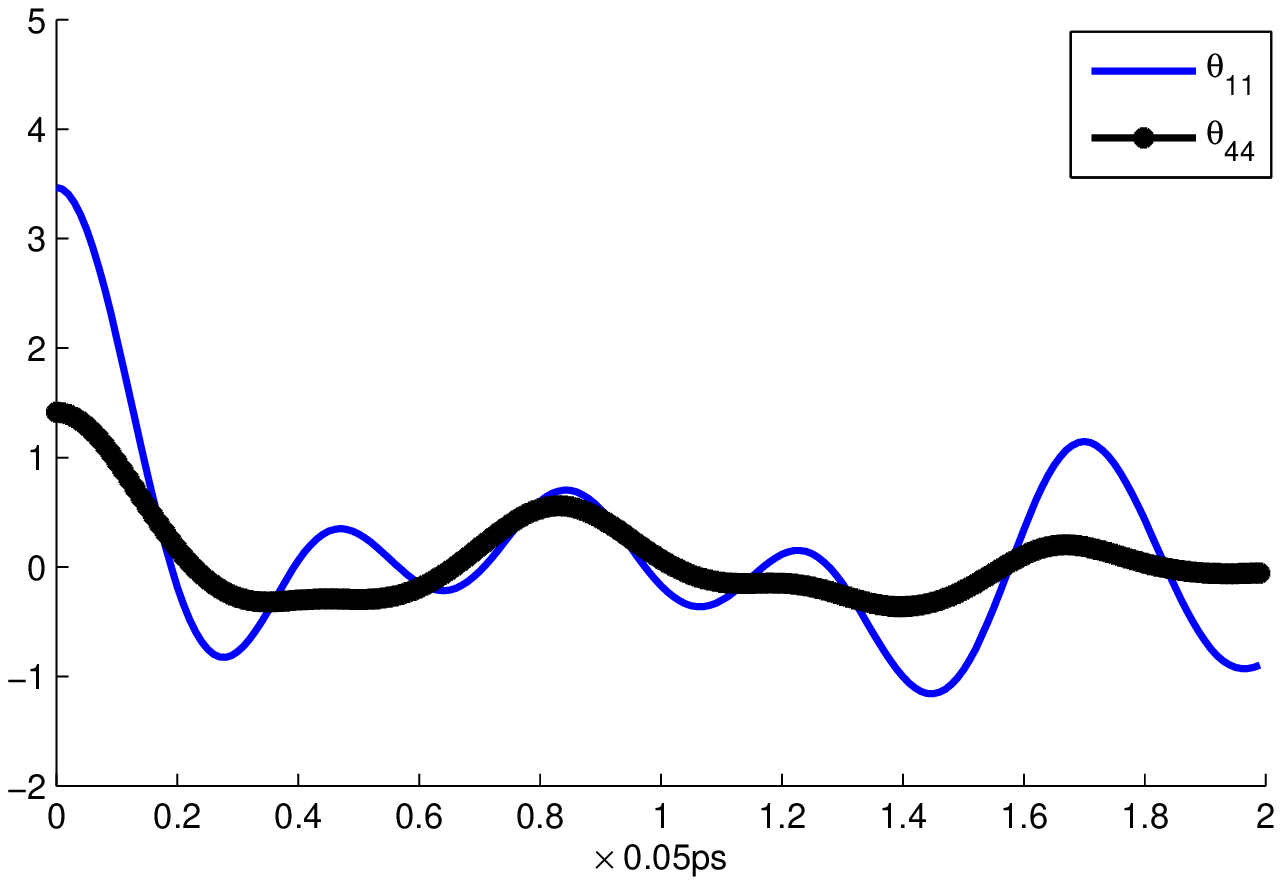}
\caption{Profiles of the two entries of the exact kernel function produced by directly computation. The coarse grained subspace are spanned by the RTB basis corresponding to the following partition. Top:  1 residue as a rigid block; Bottom: 3 residues in each block.}\label{fig7}
\end{center}
\end{figure}

Next, we consider the same type of partitions (subspace -II based on residues), but with different block sizes. In particular, we first start with a fine
partition, in which each residue is a block. We then form a coarser partition, where there are 3 residues in each block (It is clear that this partition is not based on the flexibility of the molecule). One observes from Fig. \ref{fig7} that the memory functions become smaller for the coarser partition.

\begin{figure}[htbp]
\begin{center}
\includegraphics[height=5cm, width=9cm]{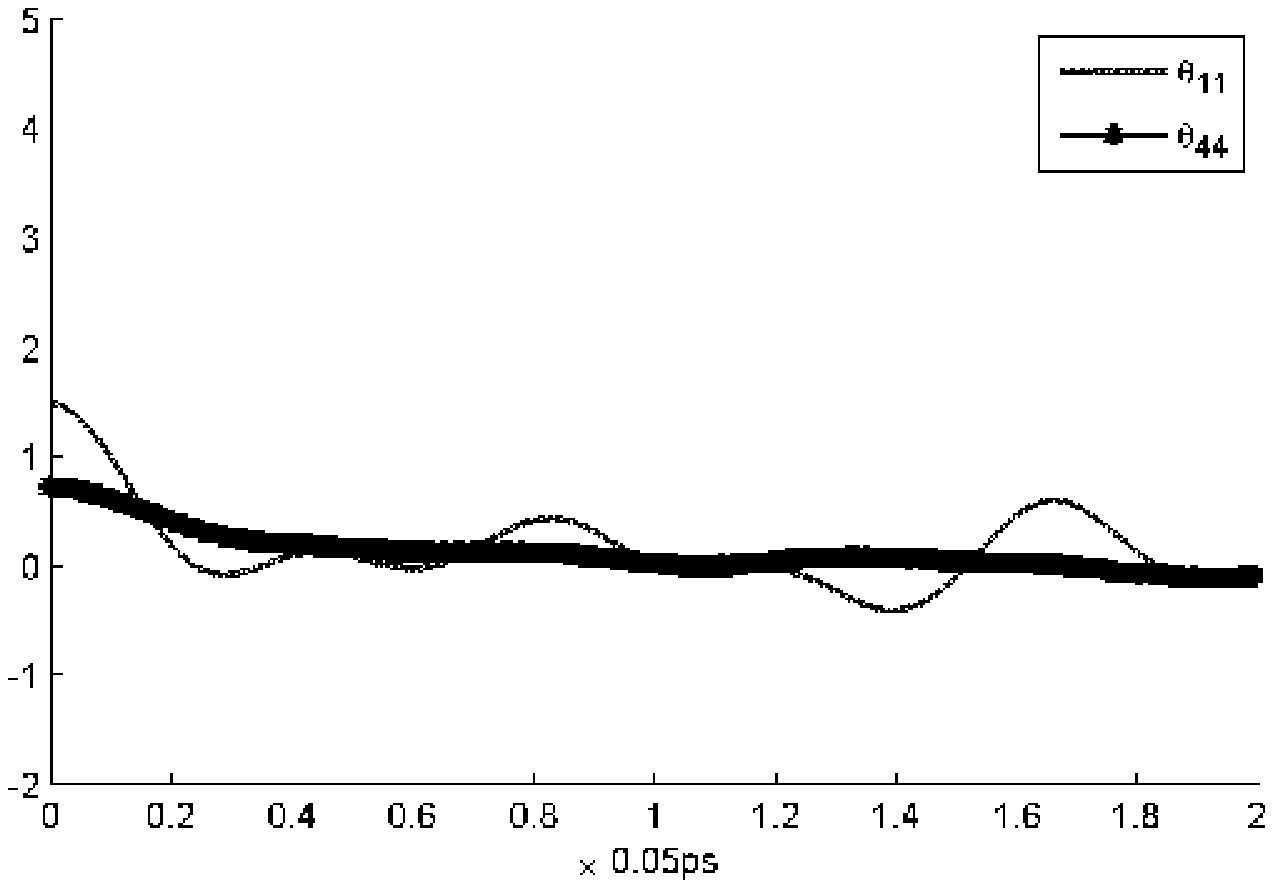}
\includegraphics[height=5cm, width=9cm]{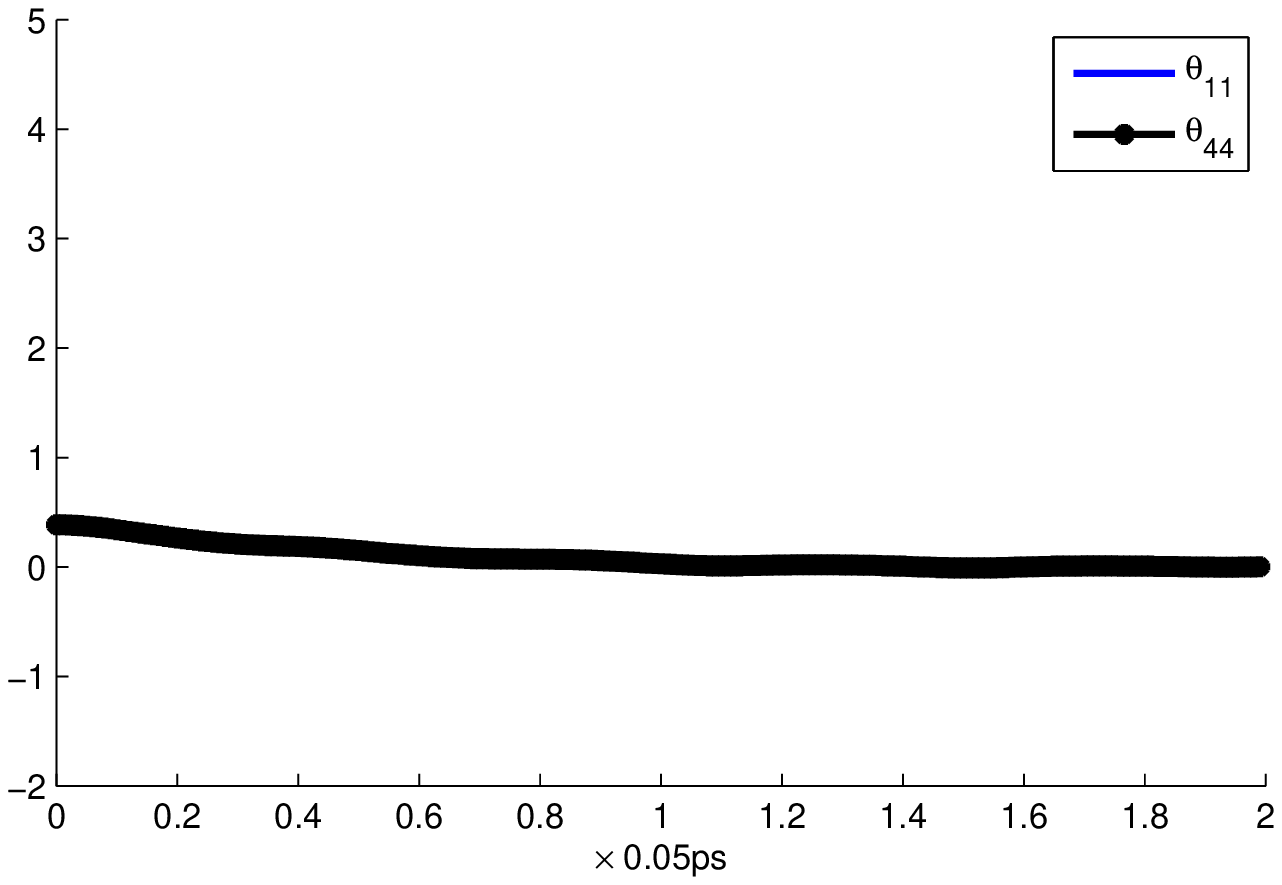}
\caption{Profiles of the first two entries of the exact kernel function produced by directly computation. The coarse grained subspace are spanned by the RTB basis corresponding to the following partition. Top 33: residues in each block; Bottom: 9 residues in one rigid block.}\label{fig8}
\end{center}
\end{figure}

To further confirm this observation, we divide the entire system equally into 22 blocks with 9 residues in each block. We also form
 a 6-block partition, each of which contains 33 residues. The results, shown in Fig. \ref{fig8}, exhibit the same trend: as we coarse-grain more  and more, the memory functions become smaller and smaller.

\section{Discussion}
In this paper, we have presented a methodology to compute memory functions which are important parameters in the generalized Langevin model. Computing such memory functions directly from molecular dynamics simulations would require
extensive effort.  In contrast, the method proposed here relies on a technique in numerical linear algebra, and it can be implemented without
performing molecular simulations.

We have also demonstrated that under the current framework, the random noise term in the generalized Langevin equation can be consistently approximated. To our knowledge, none of the existing methods offers such advantage. Together with the average force $F(\bm q)$, the generalized Langevin equation can be solved to describe the collective motion of the system. This is work in progress.

\section{Acknowledgement}
The work has been partially supported by NSF grants DMS--1109107, DMS-1216938, and DMS-1159937.
{\red This work was initialized during Chen's visitation to the Department of Mathematics, at the Pennsylvania State University.
He would like to thank the hospitality of the department. M.X. Chen was supported by the
China NSF (NSFC11301368) and the NSF of Jiangsu Province (BK20130278). }
%\bibliography{mem0}% Produces the bibliography via BibTeX.

%\nocite{*}
%\bibliography{aipsamp}% Produces the bibliography via BibTeX.

\begin{thebibliography}{48}%
\makeatletter
\providecommand \@ifxundefined [1]{%
 \@ifx{#1\undefined}
}%
\providecommand \@ifnum [1]{%
 \ifnum #1\expandafter \@firstoftwo
 \else \expandafter \@secondoftwo
 \fi
}%
\providecommand \@ifx [1]{%
 \ifx #1\expandafter \@firstoftwo
 \else \expandafter \@secondoftwo
 \fi
}%
\providecommand \natexlab [1]{#1}%
\providecommand \enquote  [1]{``#1''}%
\providecommand \bibnamefont  [1]{#1}%
\providecommand \bibfnamefont [1]{#1}%
\providecommand \citenamefont [1]{#1}%
\providecommand \href@noop [0]{\@secondoftwo}%
\providecommand \href [0]{\begingroup \@sanitize@url \@href}%
\providecommand \@href[1]{\@@startlink{#1}\@@href}%
\providecommand \@@href[1]{\endgroup#1\@@endlink}%
\providecommand \@sanitize@url [0]{\catcode `\\12\catcode `\$12\catcode
  `\&12\catcode `\#12\catcode `\^12\catcode `\_12\catcode `\%12\relax}%
\providecommand \@@startlink[1]{}%
\providecommand \@@endlink[0]{}%
\providecommand \url  [0]{\begingroup\@sanitize@url \@url }%
\providecommand \@url [1]{\endgroup\@href {#1}{\urlprefix }}%
\providecommand \urlprefix  [0]{URL }%
\providecommand \Eprint [0]{\href }%
\providecommand \doibase [0]{http://dx.doi.org/}%
\providecommand \selectlanguage [0]{\@gobble}%
\providecommand \bibinfo  [0]{\@secondoftwo}%
\providecommand \bibfield  [0]{\@secondoftwo}%
\providecommand \translation [1]{[#1]}%
\providecommand \BibitemOpen [0]{}%
\providecommand \bibitemStop [0]{}%
\providecommand \bibitemNoStop [0]{.\EOS\space}%
\providecommand \EOS [0]{\spacefactor3000\relax}%
\providecommand \BibitemShut  [1]{\csname bibitem#1\endcsname}%
\let\auto@bib@innerbib\@empty
%</preamble>
\bibitem [{\citenamefont {Schlick}(2002)}]{Schlick2002}%
  \BibitemOpen
  \bibfield  {author} {\bibinfo {author} {\bibfnamefont {T.}~\bibnamefont
  {Schlick}},\ }\href@noop {} {\emph {\bibinfo {title} {Molecular Modeling and
  Simulation: An Interdisciplinary Guide}}}\ (\bibinfo  {publisher}
  {Springer-Verlag},\ \bibinfo {year} {2002})\BibitemShut {NoStop}%
\bibitem [{\citenamefont {Leach}(2001)}]{Leach01}%
  \BibitemOpen
  \bibfield  {author} {\bibinfo {author} {\bibfnamefont {A.}~\bibnamefont
  {Leach}},\ }\href@noop {} {\emph {\bibinfo {title} {Molecular Modelling:
  Principles and Applications}}}\ (\bibinfo  {publisher} {Prentice Hall},\
  \bibinfo {year} {2001})\BibitemShut {NoStop}%
\bibitem [{\citenamefont {Baaden}\ and\ \citenamefont
  {Marrink}(2013)}]{baaden_coarse-grain_2013}%
  \BibitemOpen
  \bibfield  {author} {\bibinfo {author} {\bibfnamefont {M.}~\bibnamefont
  {Baaden}}\ and\ \bibinfo {author} {\bibfnamefont {S.~J.}\ \bibnamefont
  {Marrink}},\ }\bibfield  {title} {\enquote {\bibinfo {title} {Coarse-grain
  modelling of protein-protein interactions},}\ }\href@noop {} {\bibfield
  {journal} {\bibinfo  {journal} {Curr. Opin. Struct. Biol.}\ }\textbf
  {\bibinfo {volume} {In press.}} (\bibinfo {year} {2013})}\BibitemShut
  {NoStop}%
\bibitem [{\citenamefont {Gohlke}\ and\ \citenamefont
  {Thorpe}(2006)}]{gohlke_natural_2006}%
  \BibitemOpen
  \bibfield  {author} {\bibinfo {author} {\bibfnamefont {H.}~\bibnamefont
  {Gohlke}}\ and\ \bibinfo {author} {\bibfnamefont {M.}~\bibnamefont
  {Thorpe}},\ }\bibfield  {title} {\enquote {\bibinfo {title} {A natural coarse
  graining for simulating large biomolecular motion},}\ }\href@noop {}
  {\bibfield  {journal} {\bibinfo  {journal} {Biophys. J.}\ }\textbf {\bibinfo
  {volume} {91}},\ \bibinfo {pages} {2115--2120} (\bibinfo {year}
  {2006})}\BibitemShut {NoStop}%
\bibitem [{\citenamefont {Noid}(2013)}]{noid_perspective:_2013}%
  \BibitemOpen
  \bibfield  {author} {\bibinfo {author} {\bibfnamefont {W.~G.}\ \bibnamefont
  {Noid}},\ }\bibfield  {title} {\enquote {\bibinfo {title} {Perspective:
  coarse-grained models for biomolecular systems},}\ }\href@noop {} {\bibfield
  {journal} {\bibinfo  {journal} {J. Chem. Phys.}\ }\textbf {\bibinfo {volume}
  {139}},\ \bibinfo {pages} {090901} (\bibinfo {year} {2013})}\BibitemShut
  {NoStop}%
\bibitem [{\citenamefont {Noid}\ \emph {et~al.}(2008)\citenamefont {Noid},
  \citenamefont {Chu}, \citenamefont {Ayton}, \citenamefont {Krishna},
  \citenamefont {Izvekov}, \citenamefont {Voth}, \citenamefont {Das},\ and\
  \citenamefont {Andersen}}]{noid_multiscale_2008}%
  \BibitemOpen
  \bibfield  {author} {\bibinfo {author} {\bibfnamefont {W.~G.}\ \bibnamefont
  {Noid}}, \bibinfo {author} {\bibfnamefont {J.-W.}\ \bibnamefont {Chu}},
  \bibinfo {author} {\bibfnamefont {G.~S.}\ \bibnamefont {Ayton}}, \bibinfo
  {author} {\bibfnamefont {V.}~\bibnamefont {Krishna}}, \bibinfo {author}
  {\bibfnamefont {S.}~\bibnamefont {Izvekov}}, \bibinfo {author} {\bibfnamefont
  {G.~A.}\ \bibnamefont {Voth}}, \bibinfo {author} {\bibfnamefont
  {A.}~\bibnamefont {Das}}, \ and\ \bibinfo {author} {\bibfnamefont {H.~C.}\
  \bibnamefont {Andersen}},\ }\bibfield  {title} {\enquote {\bibinfo {title}
  {The multiscale coarse-graining method. {I.} a rigorous bridge between
  atomistic and coarse-grained models},}\ }\href@noop {} {\bibfield  {journal}
  {\bibinfo  {journal} {J. Chem. Phys.}\ }\textbf {\bibinfo {volume} {128}},\
  \bibinfo {pages} {244114} (\bibinfo {year} {2008})}\BibitemShut {NoStop}%
\bibitem [{\citenamefont {Praprotnik}, \citenamefont {Site},\ and\
  \citenamefont {Kremer}(2008)}]{praprotnik_multiscale_2008}%
  \BibitemOpen
  \bibfield  {author} {\bibinfo {author} {\bibfnamefont {M.}~\bibnamefont
  {Praprotnik}}, \bibinfo {author} {\bibfnamefont {L.~D.}\ \bibnamefont
  {Site}}, \ and\ \bibinfo {author} {\bibfnamefont {K.}~\bibnamefont
  {Kremer}},\ }\bibfield  {title} {\enquote {\bibinfo {title} {Multiscale
  simulation of soft matter: from scale bridging to adaptive resolution},}\
  }\href@noop {} {\bibfield  {journal} {\bibinfo  {journal} {Annu. Rev. Phys.
  Chem.}\ }\textbf {\bibinfo {volume} {59}},\ \bibinfo {pages} {545--571}
  (\bibinfo {year} {2008})}\BibitemShut {NoStop}%
\bibitem [{\citenamefont {Riniker}, \citenamefont {Allison},\ and\
  \citenamefont {van Gunsteren}(2012)}]{riniker_developing_2012}%
  \BibitemOpen
  \bibfield  {author} {\bibinfo {author} {\bibfnamefont {S.}~\bibnamefont
  {Riniker}}, \bibinfo {author} {\bibfnamefont {J.~R.}\ \bibnamefont
  {Allison}}, \ and\ \bibinfo {author} {\bibfnamefont {W.~F.}\ \bibnamefont
  {van Gunsteren}},\ }\bibfield  {title} {\enquote {\bibinfo {title} {On
  developing coarse-grained models for biomolecular simulation: a review},}\
  }\href@noop {} {\bibfield  {journal} {\bibinfo  {journal} {Phys. Chem. Ch.
  Ph.}\ }\textbf {\bibinfo {volume} {14}},\ \bibinfo {pages} {12423} (\bibinfo
  {year} {2012})}\BibitemShut {NoStop}%
\bibitem [{\citenamefont {Rudzinski}\ and\ \citenamefont
  {Noid}(2012)}]{rudzinski_role_2012}%
  \BibitemOpen
  \bibfield  {author} {\bibinfo {author} {\bibfnamefont {J.~F.}\ \bibnamefont
  {Rudzinski}}\ and\ \bibinfo {author} {\bibfnamefont {W.~G.}\ \bibnamefont
  {Noid}},\ }\bibfield  {title} {\enquote {\bibinfo {title} {The role of
  many-body correlations in determining potentials for coarse-grained models of
  equilibrium structure},}\ }\href@noop {} {\bibfield  {journal} {\bibinfo
  {journal} {J. Phys. Chem. B}\ }\textbf {\bibinfo {volume} {116}},\ \bibinfo
  {pages} {8621--8635} (\bibinfo {year} {2012})}\BibitemShut {NoStop}%
\bibitem [{\citenamefont {Chorin}, \citenamefont {Kast},\ and\ \citenamefont
  {Kupferman}(1998)}]{ChKaKu98}%
  \BibitemOpen
  \bibfield  {author} {\bibinfo {author} {\bibfnamefont {A.~J.}\ \bibnamefont
  {Chorin}}, \bibinfo {author} {\bibfnamefont {A.}~\bibnamefont {Kast}}, \ and\
  \bibinfo {author} {\bibfnamefont {R.}~\bibnamefont {Kupferman}},\ }\bibfield
  {title} {\enquote {\bibinfo {title} {Optimal prediction of underresolved
  dynamics},}\ }\href@noop {} {\bibfield  {journal} {\bibinfo  {journal} {Proc.
  Nat. Acad. Sci. USA}\ }\textbf {\bibinfo {volume} {96}},\ \bibinfo {pages}
  {4094 -- 4098} (\bibinfo {year} {1998})}\BibitemShut {NoStop}%
\bibitem [{\citenamefont {Chorin}, \citenamefont {Hald},\ and\ \citenamefont
  {Kupferman}(2002)}]{ChHaKu02}%
  \BibitemOpen
  \bibfield  {author} {\bibinfo {author} {\bibfnamefont {A.~J.}\ \bibnamefont
  {Chorin}}, \bibinfo {author} {\bibfnamefont {O.~H.}\ \bibnamefont {Hald}}, \
  and\ \bibinfo {author} {\bibfnamefont {R.}~\bibnamefont {Kupferman}},\
  }\bibfield  {title} {\enquote {\bibinfo {title} {Optimal prediction with
  memory},}\ }\href@noop {} {\bibfield  {journal} {\bibinfo  {journal} {Phys.
  D}\ }\textbf {\bibinfo {volume} {166}},\ \bibinfo {pages} {239--257}
  (\bibinfo {year} {2002})}\BibitemShut {NoStop}%
\bibitem [{\citenamefont {Mori}(1965)}]{Mori65}%
  \BibitemOpen
  \bibfield  {author} {\bibinfo {author} {\bibfnamefont {H.}~\bibnamefont
  {Mori}},\ }\bibfield  {title} {\enquote {\bibinfo {title} {Transport,
  collective motion, and {Brownian} motion},}\ }\href@noop {} {\bibfield
  {journal} {\bibinfo  {journal} {Prog. Theor. Phys.}\ }\textbf {\bibinfo
  {volume} {33}},\ \bibinfo {pages} {423 -- 450} (\bibinfo {year}
  {1965})}\BibitemShut {NoStop}%
\bibitem [{\citenamefont {Zwanzig}(1973)}]{Zwanzig73}%
  \BibitemOpen
  \bibfield  {author} {\bibinfo {author} {\bibfnamefont {R.}~\bibnamefont
  {Zwanzig}},\ }\bibfield  {title} {\enquote {\bibinfo {title} {Nonlinear
  generalized {Langevin} equations},}\ }\href@noop {} {\bibfield  {journal}
  {\bibinfo  {journal} {J. Stat. Phys.}\ }\textbf {\bibinfo {volume} {9}},\
  \bibinfo {pages} {215 -- 220} (\bibinfo {year} {1973})}\BibitemShut {NoStop}%
\bibitem [{\citenamefont {Berkowitz}, \citenamefont {Morgan},\ and\
  \citenamefont {McCammon}(1983)}]{berkowitz1983generalized}%
  \BibitemOpen
  \bibfield  {author} {\bibinfo {author} {\bibfnamefont {M.}~\bibnamefont
  {Berkowitz}}, \bibinfo {author} {\bibfnamefont {J.}~\bibnamefont {Morgan}}, \
  and\ \bibinfo {author} {\bibfnamefont {J.~A.}\ \bibnamefont {McCammon}},\
  }\bibfield  {title} {\enquote {\bibinfo {title} {{Generalized Langevin}
  dynamics simulations with arbitrary time-dependent memory kernels},}\
  }\href@noop {} {\bibfield  {journal} {\bibinfo  {journal} {J. Chem. Phys.}\
  }\textbf {\bibinfo {volume} {78}},\ \bibinfo {pages} {3256} (\bibinfo {year}
  {1983})}\BibitemShut {NoStop}%
\bibitem [{\citenamefont {Bu}\ and\ \citenamefont
  {Straub}(2003)}]{bu2003vibrational}%
  \BibitemOpen
  \bibfield  {author} {\bibinfo {author} {\bibfnamefont {L.}~\bibnamefont
  {Bu}}\ and\ \bibinfo {author} {\bibfnamefont {J.~E.}\ \bibnamefont
  {Straub}},\ }\bibfield  {title} {\enquote {\bibinfo {title} {Vibrational
  frequency shifts and relaxation rates for a selected vibrational mode in
  cytochrome c},}\ }\href@noop {} {\bibfield  {journal} {\bibinfo  {journal}
  {Biophys. J.}\ }\textbf {\bibinfo {volume} {85}},\ \bibinfo {pages}
  {1429--1439} (\bibinfo {year} {2003})}\BibitemShut {NoStop}%
\bibitem [{\citenamefont {Izvekov}\ and\ \citenamefont {Voth}(2006)}]{IzVo06}%
  \BibitemOpen
  \bibfield  {author} {\bibinfo {author} {\bibfnamefont {S.}~\bibnamefont
  {Izvekov}}\ and\ \bibinfo {author} {\bibfnamefont {G.~A.}\ \bibnamefont
  {Voth}},\ }\bibfield  {title} {\enquote {\bibinfo {title} {Modeling real
  dynamics in the coarse-grained representation of condensed phase systems},}\
  }\href@noop {} {\bibfield  {journal} {\bibinfo  {journal} {J. Chem. Phys.}\
  }\textbf {\bibinfo {volume} {125}},\ \bibinfo {pages} {151101--151104}
  (\bibinfo {year} {2006})}\BibitemShut {NoStop}%
\bibitem [{\citenamefont {Kamberaj}(2011)}]{kamberaj2011theoretical}%
  \BibitemOpen
  \bibfield  {author} {\bibinfo {author} {\bibfnamefont {H.}~\bibnamefont
  {Kamberaj}},\ }\bibfield  {title} {\enquote {\bibinfo {title} {A theoretical
  model for the collective motion of proteins by means of principal component
  analysis},}\ }\href@noop {} {\bibfield  {journal} {\bibinfo  {journal} {Cent.
  Eur. J. Phys.}\ }\textbf {\bibinfo {volume} {9}},\ \bibinfo {pages} {96--109}
  (\bibinfo {year} {2011})}\BibitemShut {NoStop}%
\bibitem [{\citenamefont {Lange}\ and\ \citenamefont
  {Grubm{\"u}ller}(2006)}]{lange2006collective}%
  \BibitemOpen
  \bibfield  {author} {\bibinfo {author} {\bibfnamefont {O.~F.}\ \bibnamefont
  {Lange}}\ and\ \bibinfo {author} {\bibfnamefont {H.}~\bibnamefont
  {Grubm{\"u}ller}},\ }\bibfield  {title} {\enquote {\bibinfo {title}
  {Collective {Langevin} dynamics of conformational motions in proteins},}\
  }\href@noop {} {\bibfield  {journal} {\bibinfo  {journal} {J. Chem. Phys.}\
  }\textbf {\bibinfo {volume} {124}},\ \bibinfo {pages} {214903} (\bibinfo
  {year} {2006})}\BibitemShut {NoStop}%
\bibitem [{\citenamefont {Oliva}\ \emph {et~al.}(2000)\citenamefont {Oliva},
  \citenamefont {Daura}, \citenamefont {Querol}, \citenamefont {Avil{\'e}s},\
  and\ \citenamefont {Tapia}}]{oliva2000generalized}%
  \BibitemOpen
  \bibfield  {author} {\bibinfo {author} {\bibfnamefont {B.}~\bibnamefont
  {Oliva}}, \bibinfo {author} {\bibfnamefont {X.}~\bibnamefont {Daura}},
  \bibinfo {author} {\bibfnamefont {E.}~\bibnamefont {Querol}}, \bibinfo
  {author} {\bibfnamefont {F.~X.}\ \bibnamefont {Avil{\'e}s}}, \ and\ \bibinfo
  {author} {\bibfnamefont {O.}~\bibnamefont {Tapia}},\ }\bibfield  {title}
  {\enquote {\bibinfo {title} {A generalized langevin dynamics approach to
  model solvent dynamics effects on proteins via a solvent-accessible surface.
  the carboxypeptidase a inhibitor protein as a model},}\ }\href@noop {}
  {\bibfield  {journal} {\bibinfo  {journal} {Theor. Chem. Acc.}\ }\textbf
  {\bibinfo {volume} {105}},\ \bibinfo {pages} {101--109} (\bibinfo {year}
  {2000})}\BibitemShut {NoStop}%
\bibitem [{\citenamefont {Sagnella}, \citenamefont {Straub},\ and\
  \citenamefont {Thirumalai}(2000)}]{sagnella2000time}%
  \BibitemOpen
  \bibfield  {author} {\bibinfo {author} {\bibfnamefont {D.~E.}\ \bibnamefont
  {Sagnella}}, \bibinfo {author} {\bibfnamefont {J.~E.}\ \bibnamefont
  {Straub}}, \ and\ \bibinfo {author} {\bibfnamefont {D.}~\bibnamefont
  {Thirumalai}},\ }\bibfield  {title} {\enquote {\bibinfo {title} {Time scales
  and pathways for kinetic energy relaxation in solvated proteins: Application
  to carbonmonoxy myoglobin},}\ }\href@noop {} {\bibfield  {journal} {\bibinfo
  {journal} {J. Chem. Phys.}\ }\textbf {\bibinfo {volume} {113}},\ \bibinfo
  {pages} {7702} (\bibinfo {year} {2000})}\BibitemShut {NoStop}%
\bibitem [{\citenamefont {Stepanova}(2007)}]{stepanova_dynamics_2007}%
  \BibitemOpen
  \bibfield  {author} {\bibinfo {author} {\bibfnamefont {M.}~\bibnamefont
  {Stepanova}},\ }\bibfield  {title} {\enquote {\bibinfo {title} {Dynamics of
  essential collective motions in proteins: Theory},}\ }\href@noop {}
  {\bibfield  {journal} {\bibinfo  {journal} {Phys. Rev. E}\ }\textbf {\bibinfo
  {volume} {76}},\ \bibinfo {pages} {051918} (\bibinfo {year}
  {2007})}\BibitemShut {NoStop}%
\bibitem [{\citenamefont {Zwanzig}(1961)}]{Zwanzig61}%
  \BibitemOpen
  \bibfield  {author} {\bibinfo {author} {\bibfnamefont {R.}~\bibnamefont
  {Zwanzig}},\ }\bibfield  {title} {\enquote {\bibinfo {title} {Statistical
  mechanics of irreversiblity},}\ }\href@noop {} {\bibfield  {journal}
  {\bibinfo  {journal} {Lectures in Theoretical Physics}\ }\textbf {\bibinfo
  {volume} {3}},\ \bibinfo {pages} {106--141} (\bibinfo {year}
  {1961})}\BibitemShut {NoStop}%
\bibitem [{\citenamefont {Li}(2010)}]{Li2009c}%
  \BibitemOpen
  \bibfield  {author} {\bibinfo {author} {\bibfnamefont {X.}~\bibnamefont
  {Li}},\ }\bibfield  {title} {\enquote {\bibinfo {title} {A coarse-grained
  molecular dynamics model for crystalline solids},}\ }\href@noop {} {\bibfield
   {journal} {\bibinfo  {journal} {Int. J. Numer. Meth. Engng.}\ }\textbf
  {\bibinfo {volume} {83}},\ \bibinfo {pages} {986--997} (\bibinfo {year}
  {2010})}\BibitemShut {NoStop}%
\bibitem [{\citenamefont {Chorin}\ and\ \citenamefont {Stinis}(2005)}]{ChSt05}%
  \BibitemOpen
  \bibfield  {author} {\bibinfo {author} {\bibfnamefont {A.~J.}\ \bibnamefont
  {Chorin}}\ and\ \bibinfo {author} {\bibfnamefont {P.}~\bibnamefont
  {Stinis}},\ }\bibfield  {title} {\enquote {\bibinfo {title} {Problem
  reduction, renormalization, and memory},}\ }\href@noop {} {\bibfield
  {journal} {\bibinfo  {journal} {Comm. Appl. Math. Comp. Sc.}\ }\textbf
  {\bibinfo {volume} {1}},\ \bibinfo {pages} {1--27} (\bibinfo {year}
  {2005})}\BibitemShut {NoStop}%
\bibitem [{\citenamefont {Kubo}(1966)}]{Kubo66}%
  \BibitemOpen
  \bibfield  {author} {\bibinfo {author} {\bibfnamefont {R.}~\bibnamefont
  {Kubo}},\ }\bibfield  {title} {\enquote {\bibinfo {title} {The
  fluctuation-dissipation theorem},}\ }\href@noop {} {\bibfield  {journal}
  {\bibinfo  {journal} {Rep. Prog. Phys.}\ }\textbf {\bibinfo {volume}
  {29(1)}},\ \bibinfo {pages} {255 -- 284} (\bibinfo {year}
  {1966})}\BibitemShut {NoStop}%
\bibitem [{\citenamefont {Hij{\'o}n}, \citenamefont {Serrano},\ and\
  \citenamefont {Espa{\~n}ol}(2006)}]{hijon2006markovian}%
  \BibitemOpen
  \bibfield  {author} {\bibinfo {author} {\bibfnamefont {C.}~\bibnamefont
  {Hij{\'o}n}}, \bibinfo {author} {\bibfnamefont {M.}~\bibnamefont {Serrano}},
  \ and\ \bibinfo {author} {\bibfnamefont {P.}~\bibnamefont {Espa{\~n}ol}},\
  }\bibfield  {title} {\enquote {\bibinfo {title} {Markovian approximation in a
  coarse-grained description of atomic systems},}\ }\href@noop {} {\bibfield
  {journal} {\bibinfo  {journal} {J. Chem. Phys.}\ }\textbf {\bibinfo {volume}
  {125}},\ \bibinfo {pages} {204101} (\bibinfo {year} {2006})}\BibitemShut
  {NoStop}%
\bibitem [{\citenamefont {Kauzlari{\'c}}\ \emph {et~al.}(2011)\citenamefont
  {Kauzlari{\'c}}, \citenamefont {Meier}, \citenamefont {Espa{\~n}ol},
  \citenamefont {Succi}, \citenamefont {Greiner},\ and\ \citenamefont
  {Korvink}}]{kauzlaric2011bottom}%
  \BibitemOpen
  \bibfield  {author} {\bibinfo {author} {\bibfnamefont {D.}~\bibnamefont
  {Kauzlari{\'c}}}, \bibinfo {author} {\bibfnamefont {J.~T.}\ \bibnamefont
  {Meier}}, \bibinfo {author} {\bibfnamefont {P.}~\bibnamefont {Espa{\~n}ol}},
  \bibinfo {author} {\bibfnamefont {S.}~\bibnamefont {Succi}}, \bibinfo
  {author} {\bibfnamefont {A.}~\bibnamefont {Greiner}}, \ and\ \bibinfo
  {author} {\bibfnamefont {J.~G.}\ \bibnamefont {Korvink}},\ }\bibfield
  {title} {\enquote {\bibinfo {title} {Bottom-up coarse-graining of a simple
  graphene model: The blob picture},}\ }\href@noop {} {\bibfield  {journal}
  {\bibinfo  {journal} {J. Chem. Phys.}\ }\textbf {\bibinfo {volume} {134}},\
  \bibinfo {pages} {064106--064106} (\bibinfo {year} {2011})}\BibitemShut
  {NoStop}%
\bibitem [{\citenamefont {Kauzlari{\'c}}\ \emph {et~al.}(2012)\citenamefont
  {Kauzlari{\'c}}, \citenamefont {Espa{\~n}ol}, \citenamefont {Greiner},\ and\
  \citenamefont {Succi}}]{kauzlaric2012markovian}%
  \BibitemOpen
  \bibfield  {author} {\bibinfo {author} {\bibfnamefont {D.}~\bibnamefont
  {Kauzlari{\'c}}}, \bibinfo {author} {\bibfnamefont {P.}~\bibnamefont
  {Espa{\~n}ol}}, \bibinfo {author} {\bibfnamefont {A.}~\bibnamefont
  {Greiner}}, \ and\ \bibinfo {author} {\bibfnamefont {S.}~\bibnamefont
  {Succi}},\ }\bibfield  {title} {\enquote {\bibinfo {title} {Markovian
  dissipative coarse grained molecular dynamics for a simple 2d graphene
  model},}\ }\href@noop {} {\bibfield  {journal} {\bibinfo  {journal} {The
  Journal of chemical physics}\ }\textbf {\bibinfo {volume} {137}},\ \bibinfo
  {pages} {234103} (\bibinfo {year} {2012})}\BibitemShut {NoStop}%
\bibitem [{\citenamefont {Berkowitz}\ \emph {et~al.}(1981)\citenamefont
  {Berkowitz}, \citenamefont {Morgan}, \citenamefont {Kouri},\ and\
  \citenamefont {McCammon}}]{berkowitz1981memory}%
  \BibitemOpen
  \bibfield  {author} {\bibinfo {author} {\bibfnamefont {M.}~\bibnamefont
  {Berkowitz}}, \bibinfo {author} {\bibfnamefont {J.~D.}\ \bibnamefont
  {Morgan}}, \bibinfo {author} {\bibfnamefont {D.~J.}\ \bibnamefont {Kouri}}, \
  and\ \bibinfo {author} {\bibfnamefont {J.~A.}\ \bibnamefont {McCammon}},\
  }\bibfield  {title} {\enquote {\bibinfo {title} {Memory kernels from
  molecular dynamics},}\ }\href@noop {} {\bibfield  {journal} {\bibinfo
  {journal} {J. Chem. Phys.}\ }\textbf {\bibinfo {volume} {75}},\ \bibinfo
  {pages} {2462--2463} (\bibinfo {year} {1981})}\BibitemShut {NoStop}%
\bibitem [{\citenamefont {Balescu}(1976)}]{Balescu76}%
  \BibitemOpen
  \bibfield  {author} {\bibinfo {author} {\bibfnamefont {R.}~\bibnamefont
  {Balescu}},\ }\href@noop {} {\emph {\bibinfo {title} {Equilibrium and
  Nonequilibrium Statistical Mechanics}}}\ (\bibinfo  {publisher} {John Wiley
  \& Sons},\ \bibinfo {year} {1976})\BibitemShut {NoStop}%
\bibitem [{\citenamefont {Evans}\ and\ \citenamefont {Morriss}(2008)}]{EvMo08}%
  \BibitemOpen
  \bibfield  {author} {\bibinfo {author} {\bibfnamefont {D.~J.}\ \bibnamefont
  {Evans}}\ and\ \bibinfo {author} {\bibfnamefont {G.~P.}\ \bibnamefont
  {Morriss}},\ }\href@noop {} {\emph {\bibinfo {title} {Statistical Mechanics
  of Nonequilibrium Liquids}}}\ (\bibinfo  {publisher} {ACADEMIC PRESS},\
  \bibinfo {year} {2008})\BibitemShut {NoStop}%
\bibitem [{\citenamefont {Zwanzig}(2001)}]{Zwanzigbook}%
  \BibitemOpen
  \bibfield  {author} {\bibinfo {author} {\bibfnamefont {R.}~\bibnamefont
  {Zwanzig}},\ }\href@noop {} {\emph {\bibinfo {title} {Nonequilibrium
  Statistical Mechanics}}}\ (\bibinfo  {publisher} {Oxford University Press},\
  \bibinfo {year} {2001})\BibitemShut {NoStop}%
\bibitem [{\citenamefont {Arnol'd}(1989)}]{arnol1989mathematical}%
  \BibitemOpen
  \bibfield  {author} {\bibinfo {author} {\bibfnamefont {V.~I.}\ \bibnamefont
  {Arnol'd}},\ }\href@noop {} {\emph {\bibinfo {title} {Mathematical methods of
  classical mechanics}}},\ Vol.~\bibinfo {volume} {60}\ (\bibinfo  {publisher}
  {Springer},\ \bibinfo {year} {1989})\BibitemShut {NoStop}%
\bibitem [{\citenamefont {Li}(2014)}]{Li14}%
  \BibitemOpen
  \bibfield  {author} {\bibinfo {author} {\bibfnamefont {X.}~\bibnamefont
  {Li}},\ }\bibfield  {title} {\enquote {\bibinfo {title} {Coarse-graining
  molecular dynamics models using an extended {Galerkin} projection},}\
  }\href@noop {} {\bibfield  {journal} {\bibinfo  {journal} {Int. J. Numer.
  Meth. Engng.}\ }\textbf {\bibinfo {volume} {to appear}} (\bibinfo {year}
  {2014})}\BibitemShut {NoStop}%
\bibitem [{Note1()}]{Note1}%
  \BibitemOpen
  \bibinfo {note} {This linearization is again the reference position, in which
  case the matrix is the hessian of the potential energy. This seems to be the
  only way to obtain a linear convolutional form of the memory
  function.}\BibitemShut {Stop}%
\bibitem [{\citenamefont {Yanai}, \citenamefont {Takeuchi},\ and\ \citenamefont
  {Takane}(2011)}]{yanai2011projection}%
  \BibitemOpen
  \bibfield  {author} {\bibinfo {author} {\bibfnamefont {H.}~\bibnamefont
  {Yanai}}, \bibinfo {author} {\bibfnamefont {K.}~\bibnamefont {Takeuchi}}, \
  and\ \bibinfo {author} {\bibfnamefont {Y.}~\bibnamefont {Takane}},\
  }\href@noop {} {\emph {\bibinfo {title} {Projection Matrices}}}\ (\bibinfo
  {publisher} {Springer},\ \bibinfo {year} {2011})\BibitemShut {NoStop}%
\bibitem [{Note2()}]{Note2}%
  \BibitemOpen
  \bibinfo {note} {Otherwise a nonlinear convolution term appear. Such an
  approximation is appropriate when the system is around the reference
  position. But the error would be large during structural change, in which
  case it would be necessary to keep this term.}\BibitemShut {Stop}%
\bibitem [{\citenamefont {Doob}(1944)}]{Doob44}%
  \BibitemOpen
  \bibfield  {author} {\bibinfo {author} {\bibfnamefont {J.~L.}\ \bibnamefont
  {Doob}},\ }\bibfield  {title} {\enquote {\bibinfo {title} {The elementary
  gaussian processes},}\ }\href@noop {} {\bibfield  {journal} {\bibinfo
  {journal} {Ann. Math. Stat.}\ }\textbf {\bibinfo {volume} {15}},\ \bibinfo
  {pages} {229--282} (\bibinfo {year} {1944})}\BibitemShut {NoStop}%
\bibitem [{\citenamefont {Saad}(1992)}]{saad_analysis_1992}%
  \BibitemOpen
  \bibfield  {author} {\bibinfo {author} {\bibfnamefont {Y.}~\bibnamefont
  {Saad}},\ }\bibfield  {title} {\enquote {\bibinfo {title} {Analysis of some
  {Krylov} subspace approximations to the matrix exponential operator},}\
  }\href@noop {} {\bibfield  {journal} {\bibinfo  {journal} {{SIAM} J. Numer.
  Anal.}\ }\textbf {\bibinfo {volume} {29}},\ \bibinfo {pages} {209--228}
  (\bibinfo {year} {1992})}\BibitemShut {NoStop}%
\bibitem [{\citenamefont {Diele}, \citenamefont {Moret},\ and\ \citenamefont
  {Ragni}(2009)}]{diele_error_2009}%
  \BibitemOpen
  \bibfield  {author} {\bibinfo {author} {\bibfnamefont {F.}~\bibnamefont
  {Diele}}, \bibinfo {author} {\bibfnamefont {I.}~\bibnamefont {Moret}}, \ and\
  \bibinfo {author} {\bibfnamefont {S.}~\bibnamefont {Ragni}},\ }\bibfield
  {title} {\enquote {\bibinfo {title} {Error estimates for polynomial {Krylov}
  approximations to matrix functions},}\ }\href@noop {} {\bibfield  {journal}
  {\bibinfo  {journal} {{SIAM} J. Matrix Anal. Appl.}\ }\textbf {\bibinfo
  {volume} {30}},\ \bibinfo {pages} {1546--1565} (\bibinfo {year}
  {2009})}\BibitemShut {NoStop}%
\bibitem [{\citenamefont {Saad}(2003)}]{Saad}%
  \BibitemOpen
  \bibfield  {author} {\bibinfo {author} {\bibfnamefont {Y.}~\bibnamefont
  {Saad}},\ }\href@noop {} {\emph {\bibinfo {title} {Iterative Methods for
  Sparse Linear Systems}}},\ \bibinfo {edition} {second edition}\ ed.\
  (\bibinfo  {publisher} {SIAM},\ \bibinfo {year} {2003})\BibitemShut {NoStop}%
\bibitem [{\citenamefont {Li}\ and\ \citenamefont {Cui}(2002)}]{licu02}%
  \BibitemOpen
  \bibfield  {author} {\bibinfo {author} {\bibfnamefont {G.}~\bibnamefont
  {Li}}\ and\ \bibinfo {author} {\bibfnamefont {Q.}~\bibnamefont {Cui}},\
  }\bibfield  {title} {\enquote {\bibinfo {title} {A coarse-grained normal mode
  approach for macromolecules: an efficient implementation and application to
  $ca^{2+}$-atpase},}\ }\href@noop {} {\bibfield  {journal} {\bibinfo
  {journal} {Biophys. J.}\ }\textbf {\bibinfo {volume} {83}},\ \bibinfo {pages}
  {2457--2474} (\bibinfo {year} {2002})}\BibitemShut {NoStop}%
\bibitem [{\citenamefont {F.~Tama}\ and\ \citenamefont
  {Sanejouand}(2000)}]{TaGaMaSa00}%
  \BibitemOpen
  \bibfield  {author} {\bibinfo {author} {\bibfnamefont {O.~M.}\ \bibnamefont
  {F.~Tama}, \bibfnamefont {F.~X.~Gadea}}\ and\ \bibinfo {author}
  {\bibfnamefont {Y.}~\bibnamefont {Sanejouand}},\ }\bibfield  {title}
  {\enquote {\bibinfo {title} {Building-block approach for determining
  low-frequency normal modes of macromolecules},}\ }\href@noop {} {\bibfield
  {journal} {\bibinfo  {journal} {Proteins}\ }\textbf {\bibinfo {volume}
  {41}},\ \bibinfo {pages} {1--7} (\bibinfo {year} {2000})}\BibitemShut
  {NoStop}%
\bibitem [{\citenamefont {Demerdash}\ and\ \citenamefont
  {Mitchell}(2012)}]{DeMi12}%
  \BibitemOpen
  \bibfield  {author} {\bibinfo {author} {\bibfnamefont {O.~N.~A.}\
  \bibnamefont {Demerdash}}\ and\ \bibinfo {author} {\bibfnamefont {J.~C.}\
  \bibnamefont {Mitchell}},\ }\bibfield  {title} {\enquote {\bibinfo {title}
  {Density-cluster nma: A new protein decomposition technique for
  coarse-grained normal mode analysis},}\ }\href@noop {} {\bibfield  {journal}
  {\bibinfo  {journal} {Proteins}\ }\textbf {\bibinfo {volume} {80}},\ \bibinfo
  {pages} {1766--1779} (\bibinfo {year} {2012})}\BibitemShut {NoStop}%
\bibitem [{\citenamefont {Wilson}, \citenamefont {Decius},\ and\ \citenamefont
  {Cross}(1980)}]{Wilson80}%
  \BibitemOpen
  \bibfield  {author} {\bibinfo {author} {\bibfnamefont {E.~B.}\ \bibnamefont
  {Wilson}}, \bibinfo {author} {\bibfnamefont {J.~C.}\ \bibnamefont {Decius}},
  \ and\ \bibinfo {author} {\bibfnamefont {P.~C.}\ \bibnamefont {Cross}},\
  }\href@noop {} {\emph {\bibinfo {title} {Molecular Vibrations: The Theory of
  Infrared and {Raman} Vibrational Spectra}}}\ (\bibinfo  {publisher} {Dover},\
  \bibinfo {year} {1980})\BibitemShut {NoStop}%
\bibitem [{\citenamefont {Ye}(1996)}]{ye_adaptive_1996}%
  \BibitemOpen
  \bibfield  {author} {\bibinfo {author} {\bibfnamefont {Q.}~\bibnamefont
  {Ye}},\ }\bibfield  {title} {\enquote {\bibinfo {title} {An adaptive block
  {Lanczos} algorithm},}\ }\href@noop {} {\bibfield  {journal} {\bibinfo
  {journal} {Numer. Algorithms}\ }\textbf {\bibinfo {volume} {12}},\ \bibinfo
  {pages} {97--110} (\bibinfo {year} {1996})}\BibitemShut {NoStop}%
\bibitem [{\citenamefont {Jacobs}\ \emph {et~al.}(2001)\citenamefont {Jacobs},
  \citenamefont {Rader}, \citenamefont {Kuhn},\ and\ \citenamefont
  {Thorpe}}]{JaRaKu01}%
  \BibitemOpen
  \bibfield  {author} {\bibinfo {author} {\bibfnamefont {D.~J.}\ \bibnamefont
  {Jacobs}}, \bibinfo {author} {\bibfnamefont {A.~J.}\ \bibnamefont {Rader}},
  \bibinfo {author} {\bibfnamefont {L.~A.}\ \bibnamefont {Kuhn}}, \ and\
  \bibinfo {author} {\bibfnamefont {M.~F.}\ \bibnamefont {Thorpe}},\ }\bibfield
   {title} {\enquote {\bibinfo {title} {Protein flexibility predictions using
  graph theory},}\ }\href@noop {} {\bibfield  {journal} {\bibinfo  {journal}
  {Proteins}\ }\textbf {\bibinfo {volume} {44}},\ \bibinfo {pages} {150--165}
  (\bibinfo {year} {2001})}\BibitemShut {NoStop}%
\bibitem [{\citenamefont {Tatsuoka}(1988)}]{Tatsuoka}%
  \BibitemOpen
  \bibfield  {author} {\bibinfo {author} {\bibfnamefont {M.~M.}\ \bibnamefont
  {Tatsuoka}},\ }\href@noop {} {\emph {\bibinfo {title} {Multivariate
  Analysis}}}\ (\bibinfo  {publisher} {Macmillian, New York},\ \bibinfo {year}
  {1988})\BibitemShut {NoStop}%
\end{thebibliography}

\end{document}